\documentclass[11pt]{article}

\usepackage[margin=1in]{geometry}
\usepackage{amsmath,amsfonts,amsthm,amssymb}

\usepackage{cite}
\usepackage{subcaption}
\usepackage{graphicx}
\graphicspath{{figures/}}
\usepackage{color}
\usepackage{url}

\usepackage{todonotes}

\usepackage{verbatim}

\newcommand{\bbR}{\mathbb{R}}
\newcommand{\bbN}{\mathbb{N}}
\newcommand{\Rn}{\mathbb{R}^n}
\newcommand{\bRn}{{\mathbb{R}^n}} 
\newcommand{\Rk}{\mathbb{R}^k}

\newcommand{\Grad}{\nabla}
\newcommand{\Div}{\nabla\cdot}
\newcommand{\Curl}{\nabla\times}
\newcommand{\tr}{\hbox{tr}}

\newcommand{\Ldel}{L^{\delta,\beta}}
\newcommand{\Kdel}{K^{\delta}}
\newcommand{\cdel}{c^{\delta,\beta}}
\newcommand{\mdel}{m^{\delta,\beta}}

\newcommand{\adel}{\alpha_{\delta}}
\newcommand{\T}{\mathbf{T}}

\newcommand{\note}[1]{{\color{red} #1}}

\DeclareMathOperator{\E}{E}
\DeclareMathOperator{\Res}{Res}

\newtheorem{theorem}{Theorem}

\theoremstyle{remark}

\title{Fourier spectral methods for nonlocal models}
\author{Bacim Alali and Nathan Albin\\
\
\footnotesize{Department of Mathematics, Kansas State University, Manhattan, KS}}

\begin{document}
\maketitle

\begin{abstract}
Efficient and accurate spectral solvers for nonlocal models in any spatial dimension are presented. The approach we pursue is based on the Fourier multipliers of nonlocal Laplace operators introduced in \cite{alaliAlbin2018fourier}. 
It is demonstrated that  the Fourier multipliers, and the eigenvalues in particular, can be computed accurately and efficiently. 
This is achieved through utilizing the hypergeometric representation of the Fourier multipliers in which their computation in  $n$ dimensions reduces to the computation of a 1D smooth function given in terms of $_2F_3$.  We use this representation to develop spectral techniques to solve periodic nonlocal time-dependent problems. 
For linear problems, such as the nonlocal diffusion and nonlocal wave equations, we use the diagonalizability of the nonlocal operators to produce a semi-analytic approach.
For nonlinear problems, we present a pseudo-spectral method and apply it to solve a Brusselator model with nonlocal diffusion. Accuracy and efficiency of the spectral solvers are compared against  a finite-difference solver.
\end{abstract}



\section{Introduction}
This work concerns Fourier spectral methods for nonlocal equations involving  peridynamic-type nonlocal Laplace operators. The focus is on nonlocal operators $\Ldel$, parametrized  by a spatial nonlocality parameter $\delta$ and an integral kernel exponent $\beta$, of the form \cite{alaliAlbin2018fourier}
\begin{equation}\label{eq:nonlocal_laplacian}
  \Ldel u(x) = \cdel\int_{B_\delta(x)}\frac{u(y)-u(x)}{\|y-x\|^\beta}\;dy,
\end{equation}
where $B_\delta(x)$ is a ball in $\bRn$ and the scaling constant $\cdel$ is 
defined by
\begin{eqnarray}
\label{eq:cdel-integral}
\nonumber
\cdel &:=& \left(\frac{1}{2 n}\int_{B_\delta(0)}\frac{1}{\|z\|^{\beta-2}}\;dz\right)^{-1},\\
\label{eq:cdel-explicit}
&=& \frac{2(n+2-\beta)\Gamma\left(\frac{n}{2}+1\right)}
    {\pi^{n/2}\;\delta^{n+2-\beta}}.
\end{eqnarray}

Peridynamic operators of the form \eqref{eq:nonlocal_laplacian} have been used in different applications including nonlocal diffusion, digital image correlation, and nonlocal wave phenomena, see for example \cite{bobaru2010peridynamicheat,burch2011classical,lehoucq2015novel,seleson2013interface,madenci2014peridynamic,oterkus2014peridynamic}. Nonlocal equations that involve  Laplace-type operators have been addressed in several mathematical and numerical studies including  \cite{radu2017nonlocal,du2012analysis,aksoylu2011variational, aksoylu2017application}. These nonlocal Lapalace operators are motivated by the  peridynamic theory for continuum mechanics \cite{Silling2000, silling2007peridynamic, silling2016book} and were first introduced in nonlocal vector calculus \cite{nonlocal_calc_2013}.

Spectral methods for nonlocal equations have been  developed in \cite{du2016asymptotically,du2017fast,slevinsky2018spectral}. The work in \cite{du2016asymptotically} focused on spectral approximations for a 1-D nonlocal Allen--Cahn (NAC) equation with periodic boundary conditions, and showed that these numerical methods are asymptotically compatible in the sense that they provide convergent approximations
to both nonlocal and local NAC models. As pointed out in \cite{du2016asymptotically,du2017fast,slevinsky2018spectral}, a main challenge in developing spectral methods for nonlocal equations is to accurately and efficiently compute the eigenvalues of nonlocal operators. Examples of integral kernels for which the nonlocal eigenvalues are computed explicitly are presented in \cite{du2016asymptotically}. The work  in \cite{du2017fast} provides a more general method for the accurate computations of the  eigenvalues of nonlocal operators in 1, 2, and 3 dimensions for a class of radially symmetric kernels. The method in \cite{du2017fast} is given by a hybrid algorithm that is based on reformulating the integral representations of the eigenvalues as series expansions and as solutions to ODEs. The work in \cite{slevinsky2018spectral} presents spectral algorithms for solving  nonlocal diffusion models that involve a nonlocal Laplace-Beltrami operator on the unit sphere. These algorithms are based on the diagonalizability of nonlocal diffusion operators on the sphere in the basis of spherical harmonics.

In this work, we present efficient and accurate spectral solvers for nonlocal models in $n$ dimensions. We pursue the Fourier multipliers' approach introduced in \cite{alaliAlbin2018fourier}. The Fourier multipliers of $\Ldel$ are given by the integral representation
\begin{equation}\label{eq:multiplier-cosine0}
	m(\nu) = \cdel\int_{B_\delta(0)}\frac{\cos(\nu\cdot z)-1}{\|z\|^\beta}\;dz,
\end{equation}
for $\beta<n+2$. A key observation in our approach is that this $n$-dimensional integral in \eqref{eq:multiplier-cosine0} is realized through a 1D formula for which the multipliers satisfy
\begin{equation}\label{eq:multiplier-r}
	m(r) = -r^2\;_2F_3\left(1,\frac{n+2-\beta}{2};2,\frac{n+2}{2},\frac{n+4-\beta}{2};-\frac{1}{4}r^2\delta^2\right),
\end{equation}
where $r=\|\nu\|$. Noting that the hypergeometric function $_2F_3$ is smooth in its argument, it follows from \eqref{eq:multiplier-r} that evaluating the multipliers in $\Rn$ reduces to computing this 1D smooth function. Numerical methods for computing  hypergeometric functions  such as $_2F_3$ are available in the literature, see for example \cite{NIST:DLMF}. We adopt the implementation of  $_2F_3$   that is provided by Python's {\em mpmath} library~\cite{mpmath}. 
Efficiency and accuracy study of using  {\em mpmath} to evaluate the multipliers is presented in Section~\ref{sec:efficiency}. In addition, we provide a method for the accelerated evaluation of the eigenvalues of $\Ldel$ in Section~\ref{sec:interpolation}. 

Utilizing the computations of the eigenvalues of $\Ldel$, we demonstrate spectral techniques   to solve periodic nonlocal time-dependent problems. 
For linear problems, we diagonalize the nonlocal operators in the basis of trigonometric polynomials and provide  semi-analytic solutions. We apply this in Section \ref{sec:semianalytic} to solve nonlocal diffusion and nonlocal wave equations in 2D.
For nonlinear problems, we present a pseudo-spectral method in Section~\ref{sec:pseudo-spectral} and apply it to solve  nonlocal diffusion in the Brusselator. In Section~\ref{sec:fd-vs-spectral}, we compare the spectral solver with a finite-difference solver
for the nonlocal wave equation.

\section{Fourier multipliers}\label{sec:multipliers}
For completeness of the presentation, we include a summary of the Fourier multipliers' results, developed in \cite{alaliAlbin2018fourier}, relevant to the computational spectral methods presented in Sections~\ref{sec:eigenvalues} and \ref{sec:sepctral-solvers}.
The multipliers $m$ of $\Ldel$ are defined through the Fourier transform by
\begin{equation}\label{eq:Ldel_transform}
	\Ldel u(x) = \frac{1}{(2\pi)^n}\int_{\mathbb{R}^n}m(\nu)\hat{u}(\nu)e^{i\nu\cdot x}\;d\nu,
\end{equation}
where $m$ is given by the integral representation \eqref{eq:multiplier-cosine0}. The hypergeometric representation of the multipliers is provided by the following.
\begin{theorem}\label{thm:multipliers-2F3}
Let $n\ge 1$, $\delta>0$ and $\beta<n+2$.  The Fourier multipliers  can be written in the form
\begin{equation}\label{eq:multiplier-general}
	m(\nu) = -\|\nu\|^2\,_2F_3\left(1,\frac{n+2-\beta}{2};2,\frac{n+2}{2},\frac{n+4-\beta}{2};-\frac{1}{4}\|\nu\|^2\delta^2\right).
\end{equation}
\end{theorem}
The smoothness of the $_2F_3$ function in its parameters, allows us to extend formula ~\eqref{eq:multiplier-general} for all $\beta\in\mathbb{R}\setminus\{n+4,n+6,n+8,\ldots\}$. Therefore, we utilize ~\eqref{eq:multiplier-general} to extend the definition of the multipliers $m(\nu)$ and, consequently, the definition of the operator $\Ldel$ acting on  the class of  Schwartz functions through the inverse Fourier transform in \eqref{eq:Ldel_transform}
to the case  when $\beta\ge n+2$ (with $\beta\neq n+4,n+6,n+8,\ldots$). 

A characterization of the asymptotic behavior of the multipliers, and in particular that of the  eigenvalues of $\Ldel$, is given by the following result.
\begin{theorem}\label{thm:asymptotics}
Let $n\ge 1$, $\delta>0$ and $\beta\in\mathbb{R}\setminus\{n+2,n+4,n+6,\ldots\}$.  Then, as $\|\nu\|\to\infty$,
\begin{equation*}
m(\nu) \sim 
\begin{cases}
-\frac{2n(n+2-\beta)}{\delta^2(n-\beta)}
+ 2\left(\frac{2}{\delta}\right)^{n+2-\beta}
\frac{\Gamma\left(\frac{n+4-\beta}{2}\right)\Gamma\left(\frac{n+2}{2}\right)}{(n-\beta)\Gamma\left(\frac{\beta}{2}\right)}\|\nu\|^{\beta-n}
&\text{if $\beta\ne n$},\\
-\frac{2n}{\delta^2}\left(
2\log\|\nu\|+
\log\left(\frac{\delta^2}{4}\right)+\gamma-\psi(\frac{n}{2})\right)
&\text{if $\beta = n$},
\end{cases}
\end{equation*}
where $\gamma$ is Euler's constant and $\psi$ is the digamma function.
\end{theorem}

\begin{figure}
\centering
\includegraphics[width=\textwidth]{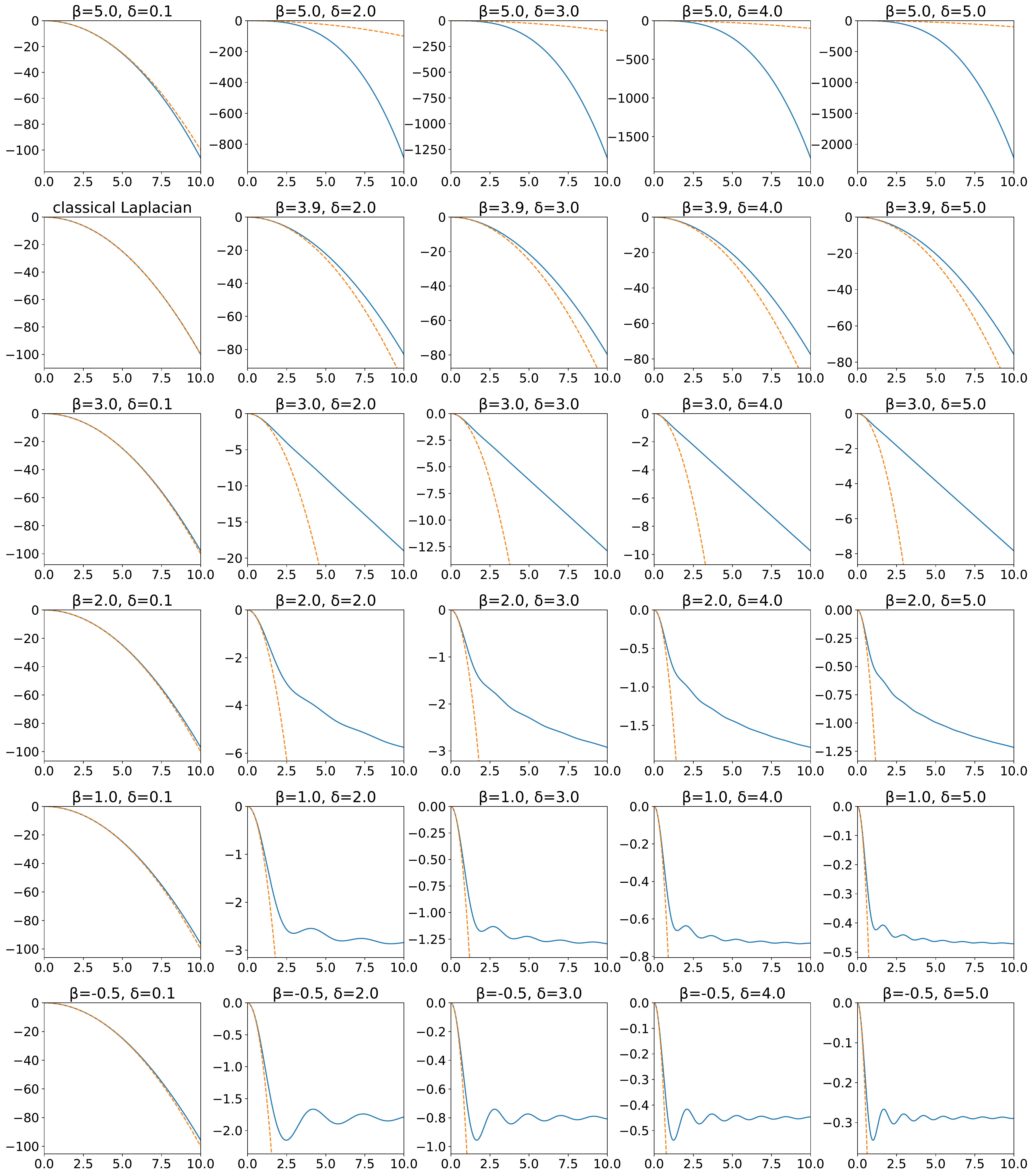}
\caption{Fourier multipliers $m(\nu)$ (vertical axis) in 2D with $\|\nu\|$ (horizontal axis) sampled at $1000$ equispaced points in the interval $[0,10]$ and with $\delta$ and $\beta$ as given in the titles.  For reference, the dashed lines show the multipliers of the classical Laplacian $m(\nu)=-\|\nu\|^2$.}
\label{fig:multiplier-comparison}
\end{figure}

Near $\nu=0$, moreover,~\eqref{eq:multiplier-general} can be used to approximate the multipliers as $m(\nu)=-\|\nu\|^2+O(\|\nu\|^4)$. Figure~\ref{fig:multiplier-comparison} provides a comparison, in the case $n=2$, of the multipliers $m(\nu)$ for various choices of $\beta$ and $\delta$.  Observe that when $\beta\approx n+2=4$ or $\delta\approx 0$, the transition to the asymptotic behavior is more gradual than when both coefficients take values away from these special cases.

\section{Eigenvalues on periodic domains}
\label{sec:eigenvalues}

The eigenvalues of $\Ldel$ can be used to develop pseudo-spectral solvers for partial integro-differential equations involving $\Ldel$ on an $n$-dimensional (flat) torus.  
In particular, suppose we consider $\Ldel$ as an operator on the periodic torus
\begin{equation*}
\T = \prod_{i=1}^n[0,\ell_i],\qquad\text{with }\ell_i>0,\quad i=1,2,\ldots,n.
\end{equation*}
For any $\alpha\in\mathbb{Z}^n$, define
\begin{equation*}
\nu_\alpha=(2\pi \alpha_1/ \ell_1,2\pi \alpha_2/ \ell_2,\ldots,2\pi \alpha_n/ \ell_n)^T\qquad\text{and}\qquad
\phi_\alpha(x) = e^{i\nu_\alpha\cdot x}.
\end{equation*}
The functions $\{\phi_\alpha\}_{\alpha\in\mathbb{Z}^n}$ form a complete set in $L^2(\T)$.  Moreover,
\begin{equation*}
  \Ldel\phi_\alpha(x) =
  \left(
  c_\delta\int_{B_\delta(0)}\frac{e^{i\nu_\alpha\cdot z}-1}{\|z\|^\beta}\;dz
  \right)\phi_\alpha(x)
  = m(\nu_\alpha)\phi_\alpha(x),
\end{equation*}
implying that $\phi_\alpha$ is an eigenfunction of $\Ldel$ with eigenvalue $m(\nu_\alpha)$.  Thus,~\eqref{eq:multiplier-general} provides an alternative to the methods established in~\cite{du2017fast} for computing the eigenvalues of $\Ldel$ on $\T$; any method for accurately evaluating the general hypergeometric function $_2F_3$ can be used to evaluate the eigenvalues in any spatial dimension $n$ and for any choices of $\beta\in \bbR\setminus\{n+4,n+6,n+8,\ldots\}$ and $\delta>0$.

\begin{table}
\centering
\begin{tabular}{|c|c|c|c|c|c|c|}
\hline
& & & & rel err \\
$n$ & $\beta$ & $\delta$ & avg time & $\|\nu\|=318\pi$\\
\hline
1 & 0.25 & 0.1 & 2.689e-04 & 1.216e-16\\
1 & 1.00 & 0.1 & 2.958e-04 & 2.401e-17\\
1 & 1.50 & 0.1 & 3.152e-04 & 8.259e-17\\
2 & 0.75 & 0.1 & 3.076e-04 & 7.624e-17\\
2 & 2.00 & 0.1 & 2.515e-04 & 8.525e-17\\
2 & 3.00 & 0.1 & 2.879e-04 & 3.949e-18\\
3 & 1.75 & 0.1 & 3.791e-04 & 2.442e-17\\
3 & 3.00 & 0.1 & 2.813e-04 & 2.960e-17\\
3 & 4.50 & 0.1 & 2.766e-04 & 3.990e-17\\
\hline
\end{tabular}
\caption{Efficiency and accuracy study, as described in Section~\ref{sec:efficiency}, using {\em mpmath}'s implementation of $_2F_3$ to compute $m(\nu)$.}
\label{tbl:efficiency}
\end{table}

\subsection{Efficiency and accuracy using Python's mpmath library}
\label{sec:efficiency}

One implementation of the $_2F_3$ hypergeometric function is provided by {\em mpmath}~\cite{mpmath}, a Python library for arbitrary precision arithmetic.  Table~\ref{tbl:efficiency} summarizes an efficiency and accuracy study using this library.
For this study, we computed, for certain choices of $n$, $\beta$ and $\delta$, the values of $m(\nu)$ with $\|\nu\|$ at $1000$ equispaced values in the interval $[1,318\pi]$.  (The interval was chosen in order to avoid relying too heavily on special values of $_2F_3$ at rational points.)  Table~\ref{tbl:efficiency} reports the results of this study.  In each row, we report the average time (in seconds) required to compute each of the $1000$ multipliers using~\eqref{eq:multiplier-general}.  We also report the relative error at $\|\nu\|=m(318\pi)$ as compared against the value of the integral in~\eqref{eq:multiplier-cosine0}, computed to 30 digits of precision using Mathematica~\cite{mathematica}.  This study was performed on a 2.7GHz laptop computer, and {\em mpmath}'s default precision (\texttt{dps=15}) was used.

\subsection{Accelerated eigenvalue approximations}\label{sec:interpolation}

\begin{figure}
\centering
\includegraphics[width=0.49\textwidth]{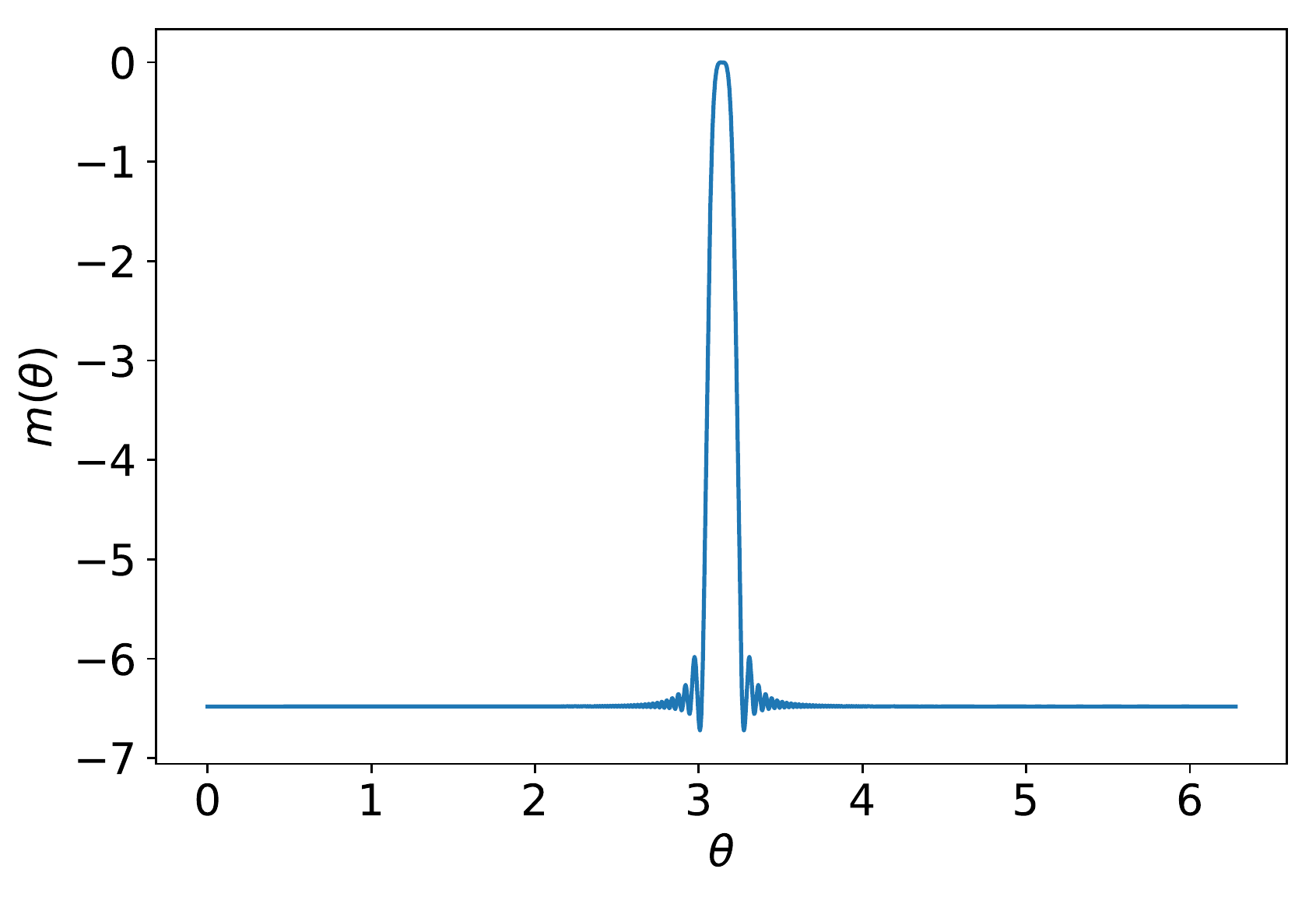}%
\includegraphics[width=0.51\textwidth]{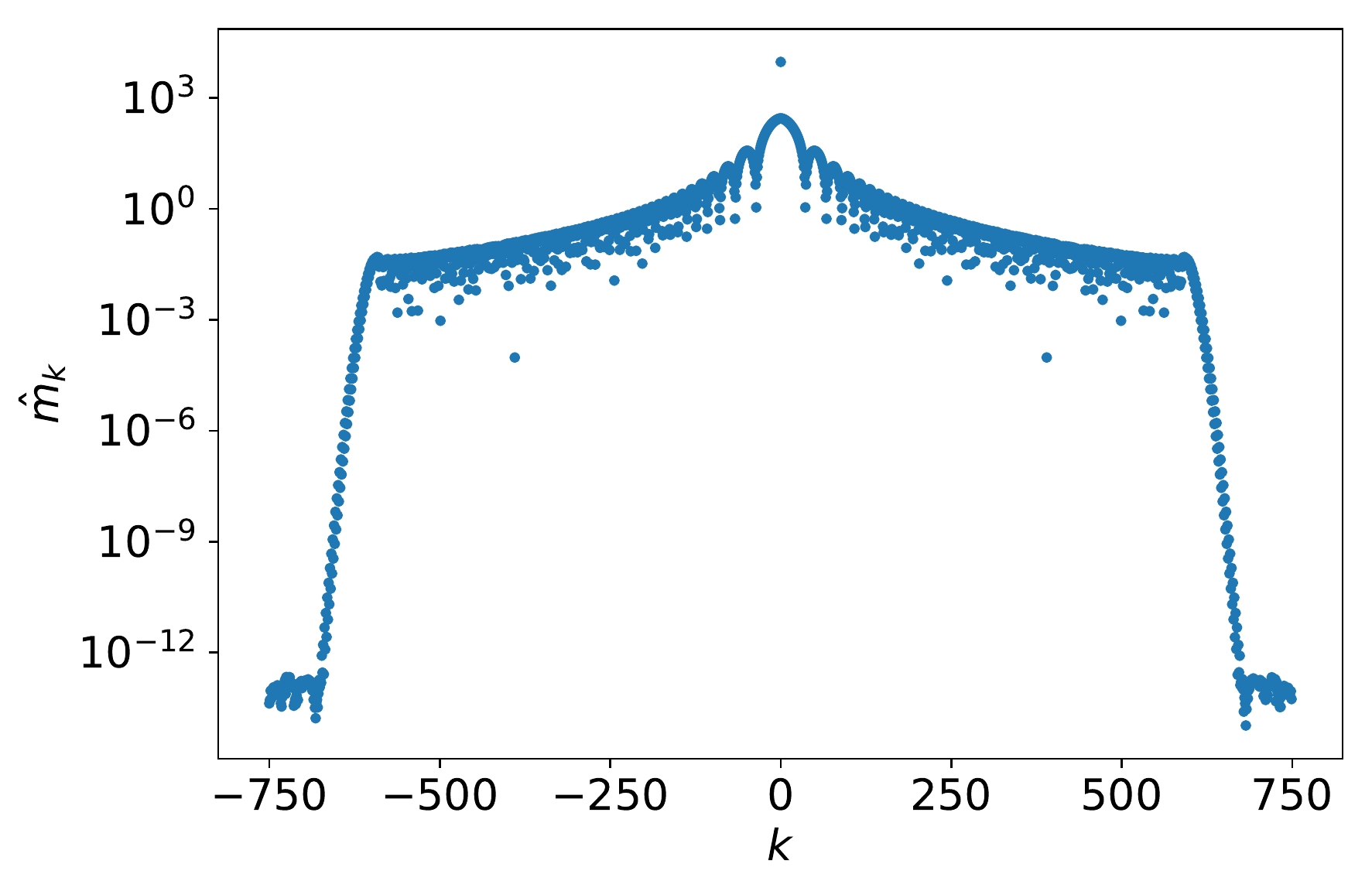}
\caption{On the left, the multipliers $m(r)$ for $n=2$, $\beta=0.5$, $\delta=1.2$ on the interval $[0,1000]$ when represented as a periodic function $m(\theta)$ as described in Section~\ref{sec:interpolation}.  On the right, the Fourier coefficients of $m(\theta)$.  $N=1500$ is sufficient to capture all significant modes.}
\label{fig:interpolation-1}
\end{figure}

\begin{figure}
\centering
\includegraphics[width=0.5\textwidth]{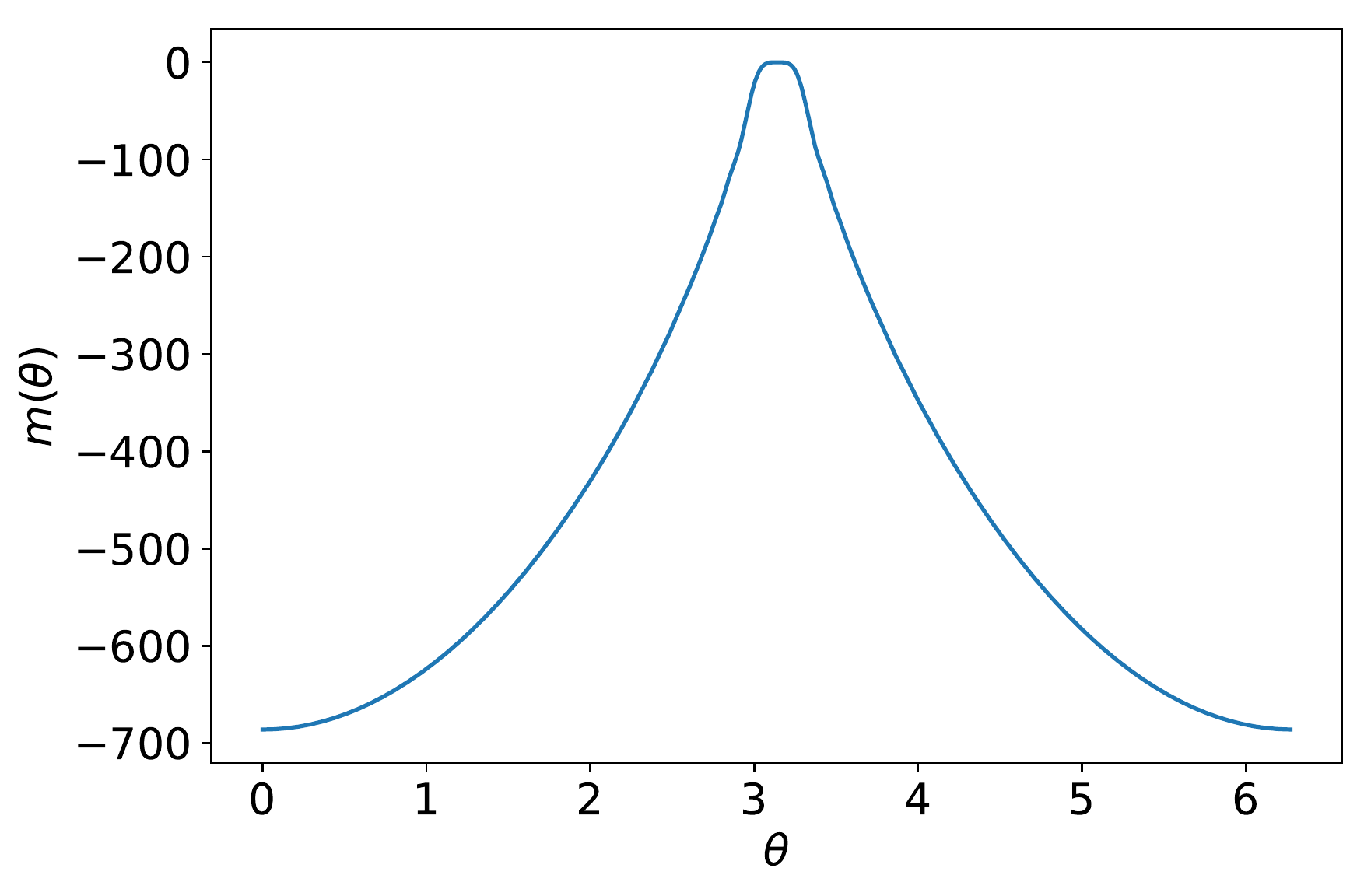}%
\includegraphics[width=0.5\textwidth]{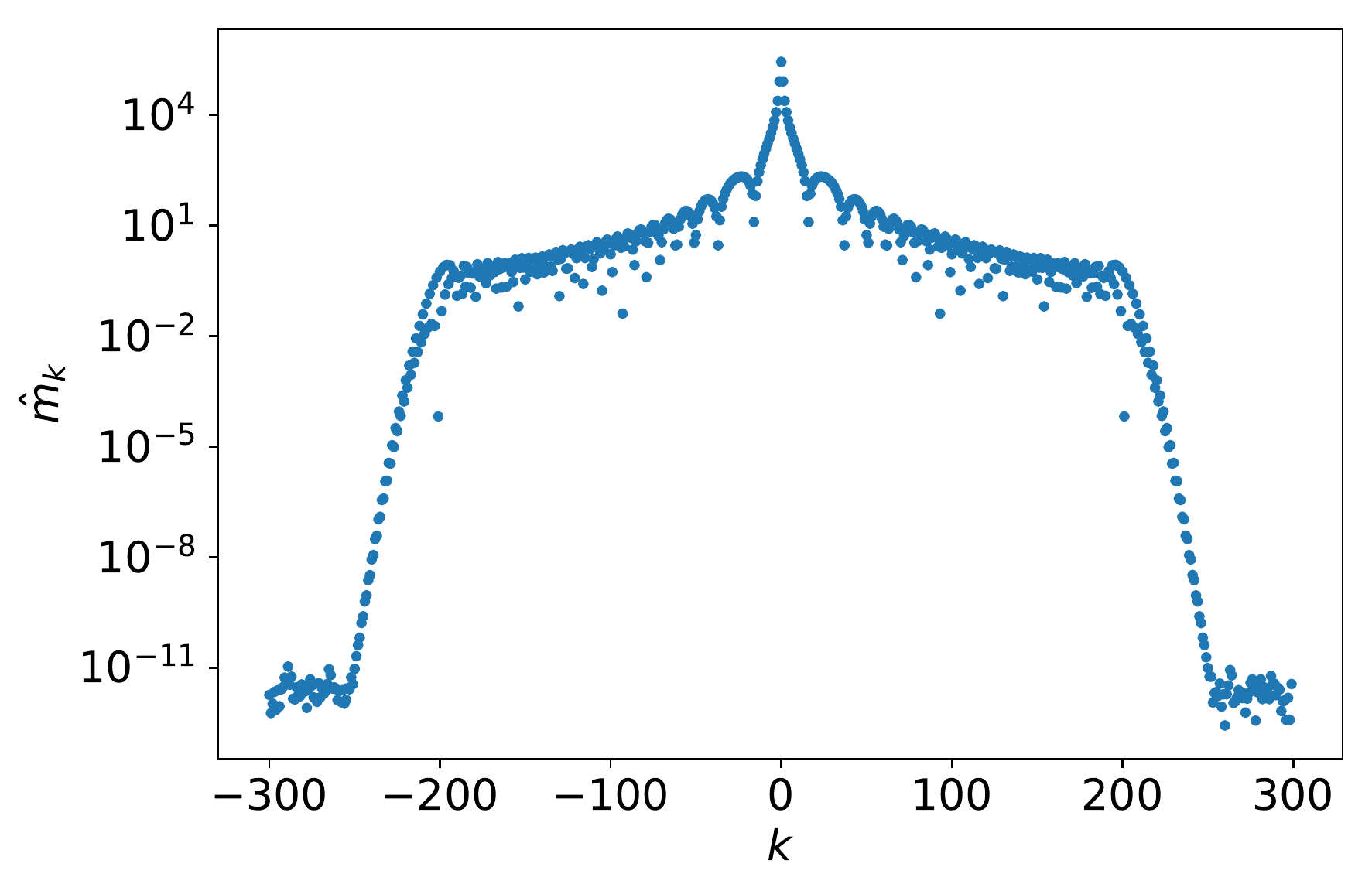}
\caption{On the left, the multipliers $m(r)$ for $n=2$, $\beta=2.3$, $\delta=0.4$ on the interval $[0,1000]$ when represented as a periodic function $m(\theta)$ as described in Section~\ref{sec:interpolation}.  On the right, the Fourier coefficients of $m(\theta)$.  $N=600$ is sufficient to capture all significant modes.}
\label{fig:interpolation-2}
\end{figure}

As demonstrated in Table~\ref{tbl:efficiency}, eigenvalues of $\Ldel$ can be computed at a rate of a few thousand per second, or around 200,000 per minute.  Considering the rapid convergence rate of Fourier-based solvers for smooth problems, this is probably sufficient for many problems.  However, for large 2D or 3D problems, or non-smooth problems that require many Fourier modes, it may be desirable to accelerate the computation.  In this section, we present one way to do this using high-accuracy Fourier interpolation combined with cubic spline interpolation.

This technique can be used effectively because $m(\nu)$ is radially symmetric, only depending on the magnitude $r=\|\nu\|$ as described in~\eqref{eq:multiplier-r}.  Suppose we wish to evaluate $m$ for many values of $r$ in an interval $[0,K]$.  Through the change of variables
\begin{equation}\label{eq:r-theta}
	r = \frac{K}{2}(1+\cos\theta),
\end{equation}
we can transform $m$ into a periodic function of $\theta$.  Moreover, since $m(r)$ is an analytic function, so is $m(\theta)$, so $m(\theta)$ has a Fourier series
\begin{equation*}
	m(\theta) = \sum_{k=-\infty}^\infty \hat{m}_ke^{ik\theta}
\end{equation*}
with rapidly decaying Fourier coefficients $\hat{m}_k$.  Figures~\ref{fig:interpolation-1} and~\ref{fig:interpolation-2} show two examples.

The technique for accelerating the evaluation of $m(r)$ proceeds as follows.  First, we sample $m(\theta)$ at the points $\theta=2\pi j/N$ for $j=0,1,\ldots,N-1$.
The parameter $N$ should be chosen sufficiently large that the coefficients $\hat{m}_k$ are very small for $|k|>N/2$.  Now let $M>N$ be an even integer.  $m(\theta)$ can be resampled with high accuracy on the points $\theta=2\pi j/M$ for $j=0,1,\ldots,M-1$ using zero padding.  (One simply uses an FFT on the $N$ sample points, pads the $N$ coefficients with $M-N$ zeros in the highest frequency modes, and transforms back with an inverse FFT.)  Because of the efficiency of the FFT, this can be done very quickly.  Finally, by inverting the transformation~\eqref{eq:r-theta}, we obtain a highly accurate approximation of $m(z)$ at the points $r=\frac{K}{2}(1+\cos(2\pi j/M))$ for $j=0,1,\ldots,M/2$, the Chebyshev nodes of the second kind in the interval $[0,K]$.  If $M$ is chosen sufficiently large, the cubic spline interpolation of these points will provide an accurate approximation of $m(r)$ throughout the interval.

\begin{table}
\centering
\begin{tabular}{|l|r|r|}
	\hline
	$n$ & 2 & 2 \\
    $\beta$ & 0.5 & 2.3 \\
    $\delta$ & 1.2 & 0.4 \\
    $K$ & 1000 & 1000 \\
    $N$ & 1500 & 600 \\
    $M$ & 20000 & 10000 \\
    prep time & 3.117 s & 0.423 s \\
    avg interp time & 8.720e-8 s & 1.354e-7 s \\
    avg direct time & 1.949e-3 s & 8.099e-4 s \\
    max error & 3.758e-10 & 5.301e-11 \\
    \hline
\end{tabular}
\caption{Results of the interpolation test described in Section~\ref{sec:interpolation}.}
\label{tab:interpolation}
\end{table}

Table~\ref{tab:interpolation} demonstrates the effectiveness of this approach for the two cases shown in Figures~\ref{fig:interpolation-1} and~\ref{fig:interpolation-2}.  Along with the relevant coefficients, the ``prep time'' row indicates the time required to build the interpolation data (including sampling the function on the $N$-point grid, the FFT-based interpolation from the $N$-point grid to the $M$-point grid, and the construction of the cubic spline interpolant on the $M$-point grid), the ``avg interp time'' row reports the average time required to evaluate $m(\nu)$ using the cubic interpolation after the interpolant is constructed, and the ``avg direct time'' row reports the average time required to evaluate $m(\nu)$ using~\eqref{eq:multiplier-general} directly.  The row titled ``max error'' reports the maximum value of $|x_t-x_a|/(1+|x_t|)$ where $x_a$ is the interpolated value and $x_t$ is the value from the formula in~\eqref{eq:multiplier-general} at 10,000 uniformly spaced points in $[0,K]$.  Again, this study was performed on a 2.7GHz laptop computer using {\em mpmath} to compute the values of $_2F_3$.  The FFT and cubic interpolation were provided by {\em scipy}~\cite{scipy}.  As can be seen from this example, by spending a few seconds generating interpolation data, it is possible to increase the speed of computing the multipliers by several orders of magnitude.

\section{Spectral methods for time-dependent periodic problems}
\label{sec:sepctral-solvers}
In this section, we demonstrate techniques for using the eigenvalues of $\Ldel$ to solve time-dependent problems with nonlocal spatial operators.  We also compare these solvers with one of the finite-difference-like solvers described in~\cite{tian2013analysis}. 

\subsection{Semi-analytic solutions to linear equations}\label{sec:semianalytic}

On a periodic torus $\T$, the peridynamic heat and wave equations take the form,
\begin{equation}\label{eq:nonlocal-heat-per}
	\begin{cases}
	u_t = \Ldel u\qquad\text{in }(0,\infty)\times \T,\\
    u(0,x) = u_0(x),
	\end{cases}
\end{equation}
and
\begin{equation}\label{eq:nonlocal-wave-per}
	\begin{cases}
	u_{tt} = \Ldel u\qquad\text{in }(0,\infty)\times \T,\\
    u(0,x) = u_0(x),\\
    u_t(0,x) = v_0(x),
	\end{cases} 
\end{equation}
respectively.  The solution to~\eqref{eq:nonlocal-heat-per} can be written as
\begin{equation}\label{eq:nonlocal-heat-per-sol}
	u(t,x) = \sum_{\alpha\in\mathbb{Z}^n}e^{\lambda_\alpha t}\hat{u}_{0,\alpha} e^{i\nu_\alpha\cdot x},
\end{equation}
where $\lambda_\alpha=m(\nu_\alpha)$ and $\hat{u}_{0,\alpha}$ are the Fourier coefficients for $u_0$.  Similarly, the solution to~\eqref{eq:nonlocal-wave-per} can be written as
\begin{equation}\label{eq:nonlocal-wave-per-sol}
\begin{split}
	u(t,x) &= \sum_{\alpha\in\mathbb{Z}^n}\hat{u}_\alpha(t)e^{i\nu_\alpha\cdot x},
    \qquad\text{where}\\
    \hat{u}_\alpha(t) &= 
    \frac{1}{2}\left(
    	\hat{u}_{0,\alpha}+\frac{\hat{v}_{0,\alpha}}{\sqrt{\lambda_\alpha}}
    \right)e^{\sqrt{\lambda_\alpha}t}
    +
	\frac{1}{2}\left(
    	\hat{u}_{0,\alpha}-\frac{\hat{v}_{0,\alpha}}{\sqrt{\lambda_\alpha}}
    \right)e^{-\sqrt{\lambda_\alpha}t},
\end{split}
\end{equation}
where $\lambda_\alpha=m(\nu_\alpha)$ (as before), and $\hat{u}_{0,\alpha}$ and $\hat{v}_{0,\alpha}$ are the Fourier coefficients for $u_0$ and $v_0$ respectively.  Linear nonlocal problems like these can be solved semi-analytically on $\T$ using the Fast Fourier Transform (FFT) to approximate the Fourier coefficients and the initial data.  To do so, we simply apply the FFT to the initial data, apply one of the solution formulas above, and then use the inverse FFT to find the solution at the specified time.

\subsubsection{Example: the nonlocal heat equation}\label{sec:nonlocal-heat-example}

\begin{figure}
\centering
\includegraphics[width=1.03\textwidth]{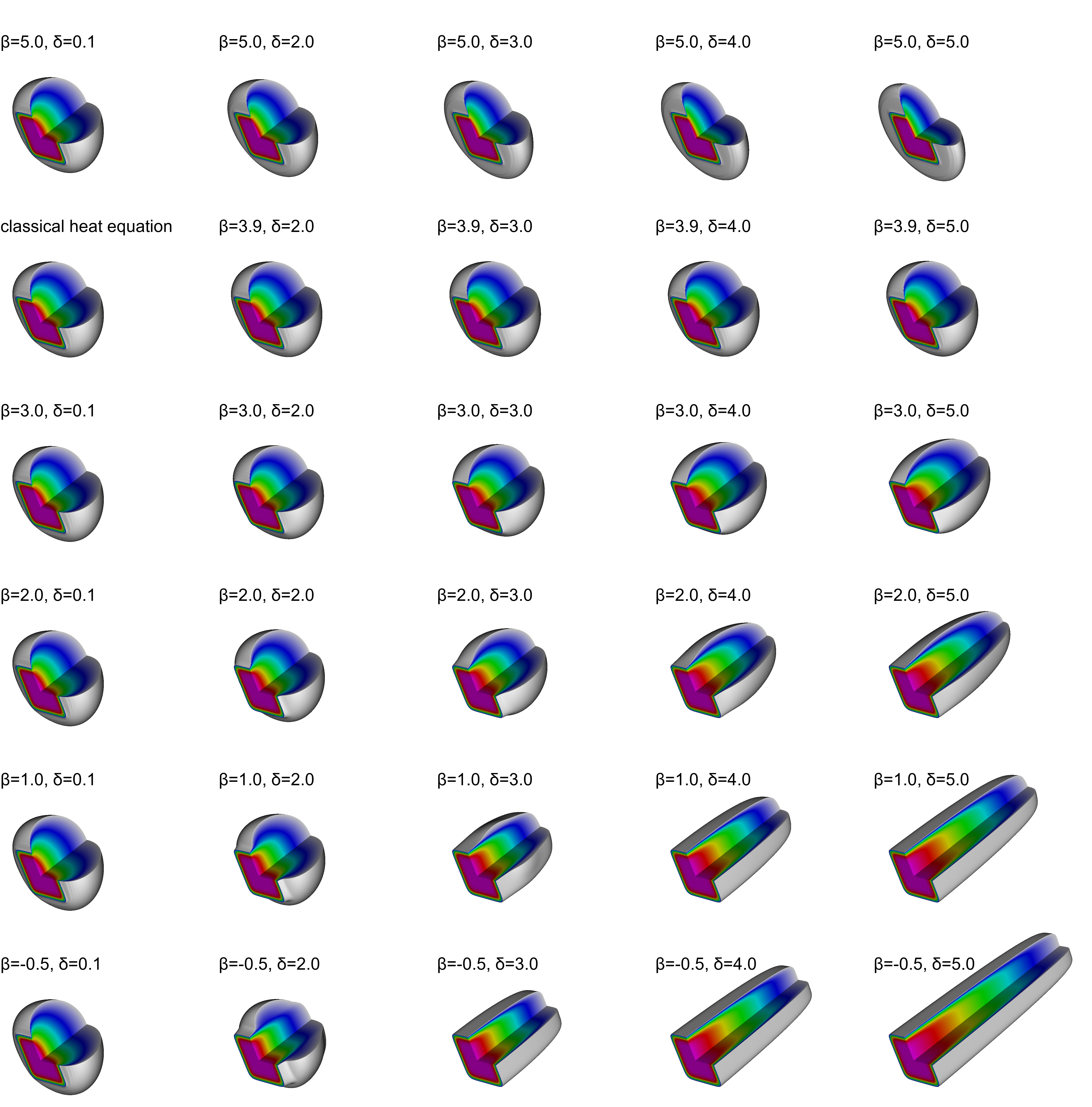}
\caption{Solutions to the nonlocal heat equation described in Section~\ref{sec:nonlocal-heat-example} with parameters $\beta$ and $\delta$ as indicated above each plot.}
\label{fig:heat-example-1}
\end{figure}

Consider, for example, the nonlocal heat equation in~\eqref{eq:nonlocal-heat-per} on the 2D torus $\T=[-10,10]^2$ with initial data $u_0(x,y)=\exp(-x^8-y^8)$.  Although this is not truly a periodic function, it is indistinguishable from a smooth periodic function at machine precision.  Figure~\ref{fig:heat-example-1} shows the solutions in $(t,x,y)$-space (with the positive $t$ direction receding back into the page and toward the upper right) as computed by the methods described in this section.  In all cases, the spatial domain was discretized using a grid with $800\times 800$ points, which is sufficient to capture the Fourier spectrum of $u_0$ to machine precision.  The Fourier coefficients for $u_0$ were approximated by FFT.  The temporal domain $[0,15]$ was discretized using 250 equally spaced points.  For each choice of operator and each time level,~\eqref{eq:nonlocal-heat-per-sol} was used to compute the solution $u(t,x,y)$ in the frequency domain and transformed back to the spatial domain by inverse FFT.  To aid in visualization, values of $u < 0.1$ are omitted from the 3D plots; the outer ``shell'' of each plot corresponds to the level surface $u(t,x,y)=0.1$.  A wedge running parallel to the $t$ axis has been removed in order to show values of $u$ near $(x,y)=(0,0)$.

Note that, as could be inferred from Figure~\ref{fig:multiplier-comparison}, when $\delta\approx 0$, the close agreement between $m(\nu)$ and the multipliers of the Laplacian ensures that the smoothed square bump decays similarly by rapidly transitioning into a circular bump and dissipating away. On the other hand, when $\delta$ is relatively large, the dependence on $\beta$ becomes more pronounced. For $\beta<n+2=4$, the solution values spread more slowly than the classical heat equation. In fact, for $\beta<n=2$, the solution does not appear to spread at all.  On the other hand, when $\beta>n+2$, the solution spreads and decays more rapidly than in the classical heat equation.

\begin{figure}
\centering
\hspace*{-.8cm}
\includegraphics[width=1.07\textwidth]{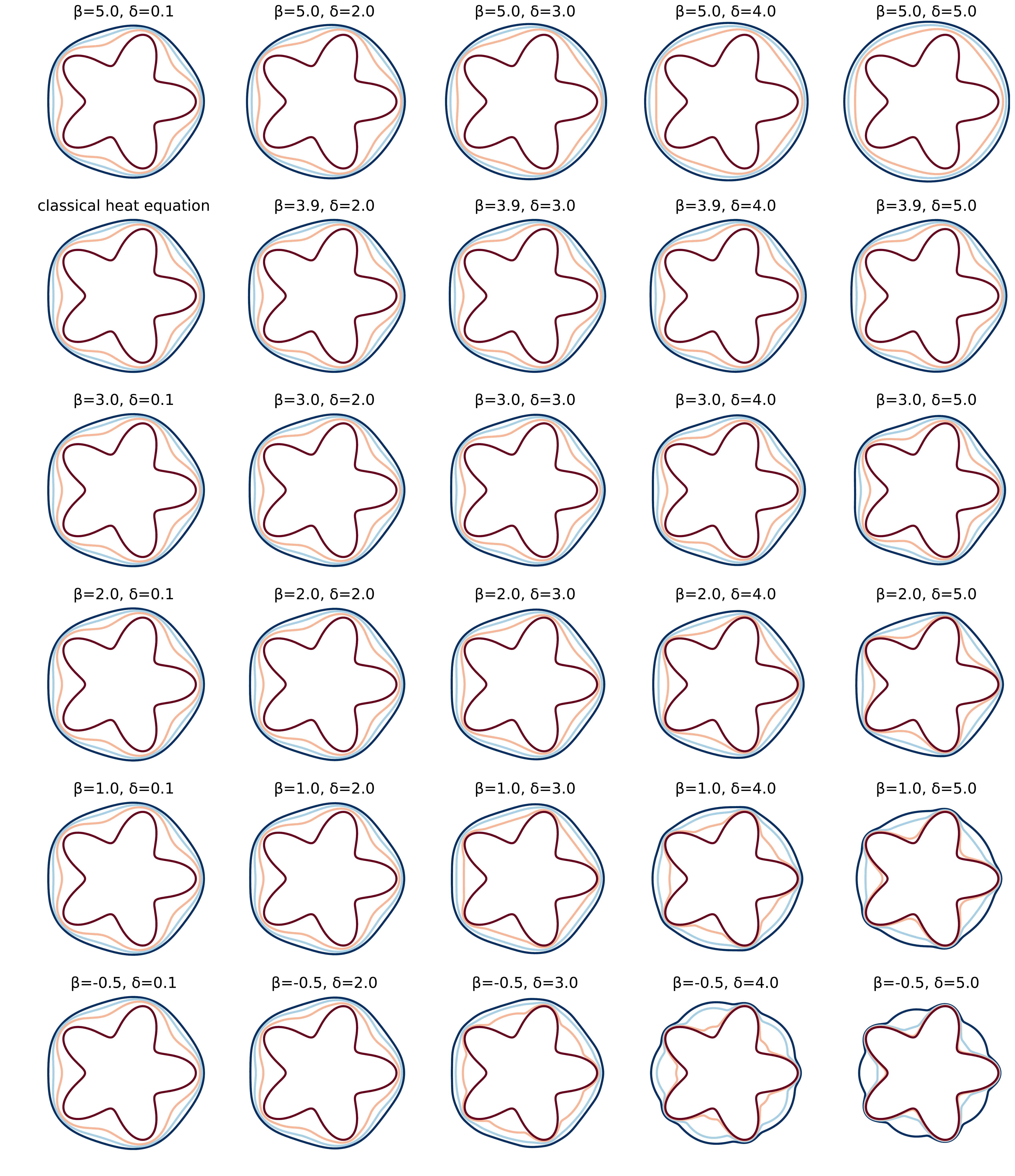}
\caption{Solutions to the nonlocal heat equation in 2D as described in Section~\ref{sec:nonlocal-heat-example} with parameters $\beta$ and $\delta$ as indicated above each plot.  The curves indicate the $u=0.1$ contour at 4 equally spaced times between $t=0$ and $t=2$.}
\label{fig:heat-example-2}
\end{figure}

A similar phenomenon can be seen in Figure~\ref{fig:heat-example-2}.  This example is identical to the previous except in the initial condition, which can be written in polar form as
\begin{equation*}
    u(t,r,\theta) = \exp\left[-\left(\frac{r}{4(1 + 0.3\sin(5\theta))}\right)^8\right].
\end{equation*}
The figure shows how the level lines $u=0.1$ change as a function of time.  Again, the results for $\delta\approx 0$ are quite similar to the results from the classical heat equation.  When $\delta$ is larger, though,  one again sees how the choice of parameter $\beta$ affects the diffusion phenomenon as compared to  classical diffusion: values of $\beta<n+2=4$ slow the diffusion while values of $\beta>n+2=4$ accelerate the diffusion.

\subsubsection{Example: the nonlocal wave equation}
\label{sec:nonlocal-wave-example}

\begin{figure}
\centering
\includegraphics[width=\textwidth]{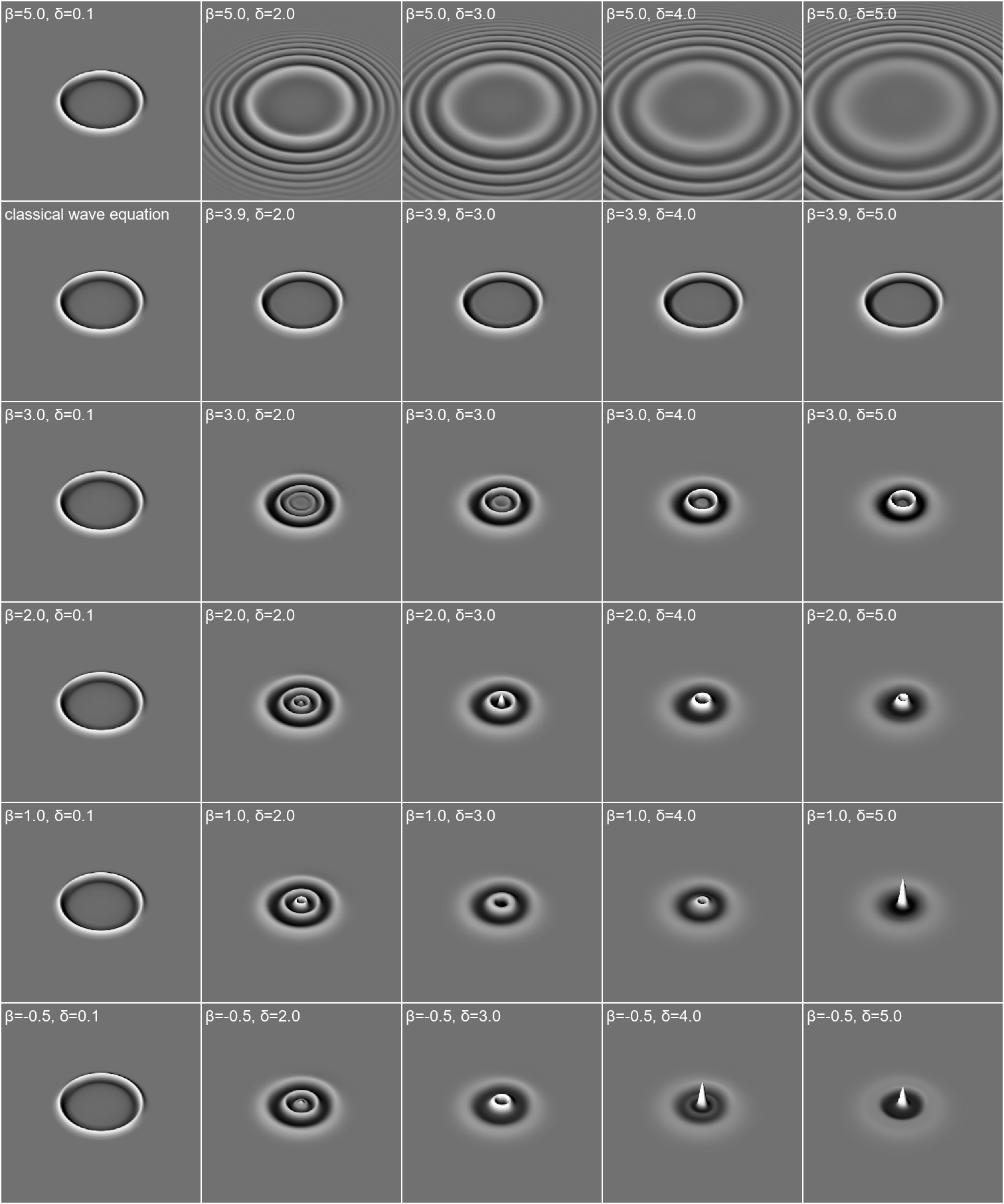}
\caption{Solutions to the nonlocal wave equation described in Section~\ref{sec:nonlocal-wave-example} with parameters $\beta$ and $\delta$ as indicated in the corner of each plot.  All plots show the solution $u$ at time $t=10$.}
\label{fig:2d-wave-example}
\end{figure}

As a second example, consider the 2D nonlocal wave equation~\eqref{eq:nonlocal-wave-per} on the torus $\T=[-48,48]^2$ with Gaussian initial conditions $u_0(x)=\exp(-\|x\|^2)$, $v_0(x)=0$.  Once again, this function is periodic to machine precision.  Figure~\ref{fig:2d-wave-example} shows the solutions $u(t,x)$ at time $t=10$ for several choices of nonlocal operator (as well as the classical wave equation).  In all cases, the spatial domain was discretized using 400 points in each direction.  (The Gaussian can be resolved using fewer points, but the final images are not as clear.)  The solution at time $t=10$ was then produced using~\eqref{eq:nonlocal-wave-per-sol} and the FFT.  As in Section~\ref{sec:nonlocal-heat-example}, we observe that when $\delta\approx 0$, the behavior is quite similar to the classical wave equation while the cases with larger $\delta$ exhibit more sensitivity to the parameter $\beta$.  Again, $\beta=n+2=4$ appears as a threshold with $\beta<4$ displaying slower wave propagation than the classical case and $\beta>4$ displaying faster.

\subsection{Pseudo-spectral methods}
\label{sec:pseudo-spectral}

For nonlinear problems (or even problems with varying coefficients), the semi-analytic method from the previous section tend not to work well.  The main difficulty arises from the fact that products of functions, under Fourier transform, give rise to convolutions that prevent the decoupling of the Fourier modes present in the linear case.

Nonetheless, as shown in~\cite{du2017fast}, the fact that the trigonometric polynomials are eigenfunctions for $\Ldel$ in the periodic setting
allows the development of fast and accurate pseudo-spectral (or Fourier collocation) numerical
methods.  Essentially, whenever we need to apply $\Ldel$ to a
function, we first pass to the frequency domain using an FFT, then
apply the operator there (which amounts to simply scaling each
Fourier mode by the appropriate eigenvalue) and then invert the FFT to
return to the spatial domain.  Thus $\Ldel$ is applied efficiently in the frequency domain, while function multiplication can be applied efficiently in the spatial domain.

\subsubsection{Nonlocal diffusion in the Brusselator}\label{sec:brusselator}

As an example, consider the following  nonlocal Brusselator reaction-diffusion model
\begin{equation*}
\begin{cases}
    u_t = D_u\Ldel u + a - (b+1)u + u^2 v,\\
    v_t = D_v\Ldel v + bu - u^2 v,
\end{cases}
\end{equation*}
where we have simply replaced the Laplacian by its nonlocal version.

\begin{figure}
\centering
\includegraphics[width=\textwidth]{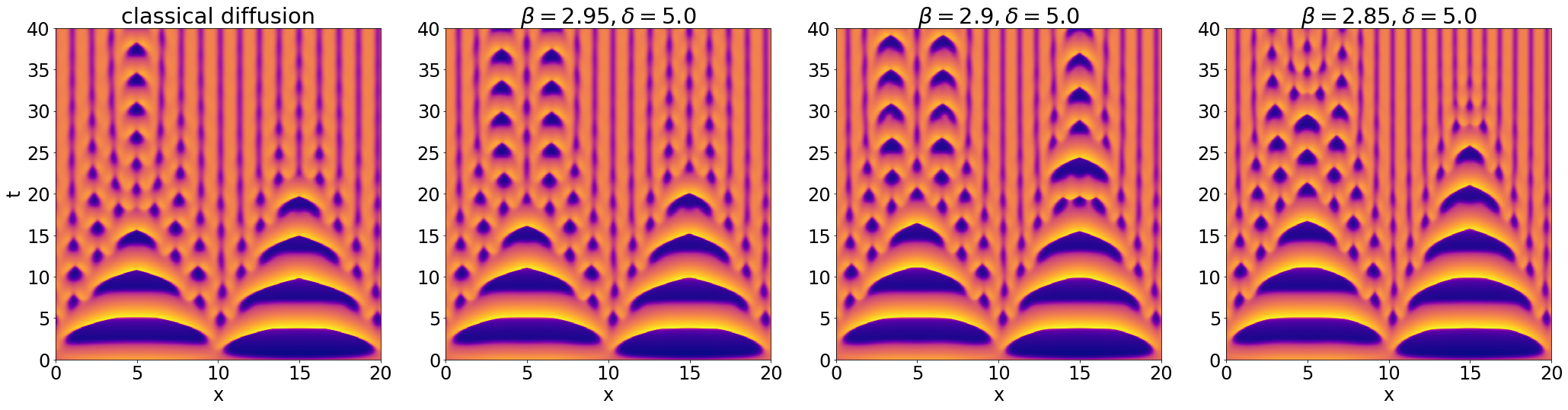}
\caption{Solutions to the Brusselator equations using classical diffusion as well as $\Ldel$ for several choices of parameters as described in Section~\ref{sec:brusselator}.  The colors in the plots indicate the values of the concentration $u$ as a function of $x$ and $t$.}
\label{fig:brusselator-example}
\end{figure}

Figure~\ref{fig:brusselator-example} shows several solutions to the 1D periodic system approximated using the pseudospectral method described in Section~\ref{sec:pseudo-spectral}. For this example, we chose parameter values
\begin{equation*}
    D_u=0.0625,\quad
    D_v=0.12,\quad
    a=3,\quad
    \text{and}\quad b=11.
\end{equation*}
These equations were solved on the periodic spatial interval $[0,20]$ for times in the interval $[0,40]$. The intial conditions chosen were
\begin{align*}
    u(0,x) &= a\left(1+\frac{1}{2}\sin\left(\frac{\pi}{10}x\right)\right) \\
    v(0,x) &=
    \frac{b}{a} + \frac{1}{10}
    \cos\left(\frac{3\pi}{5}x\right).
\end{align*}

The numerical solution was advanced using a standard RK4 timestepping scheme combined with a spectral filter (zeroing the coefficients of the 1/3 largest wave numbers) to prevent aliasing errors.  A CFL condition of the form $\Delta t=1.9\Delta x^2$ was chosen empirically from numerical tests.  The periodic interval $[0,20]$ was discretized using $1600$ points, resulting in a temporal discretization using $134,737$ points.  This discretization is sufficiently fine that the plots in the figure are not altered if the number of spatial nodes is halved.

The complexity of these figures illustrates again the importance of utilizing highly accurate methods as well as the nontrivial role played by $\beta$ in the nonlocal operator.  Small differences in this parameter can yield significant differences in solutions.

\subsection{Comparison with a finite difference method}
\label{sec:fd-vs-spectral}
\begin{figure}
  \centering
  \begin{subfigure}[t]{0.322\textwidth}
    \includegraphics[width=\textwidth]{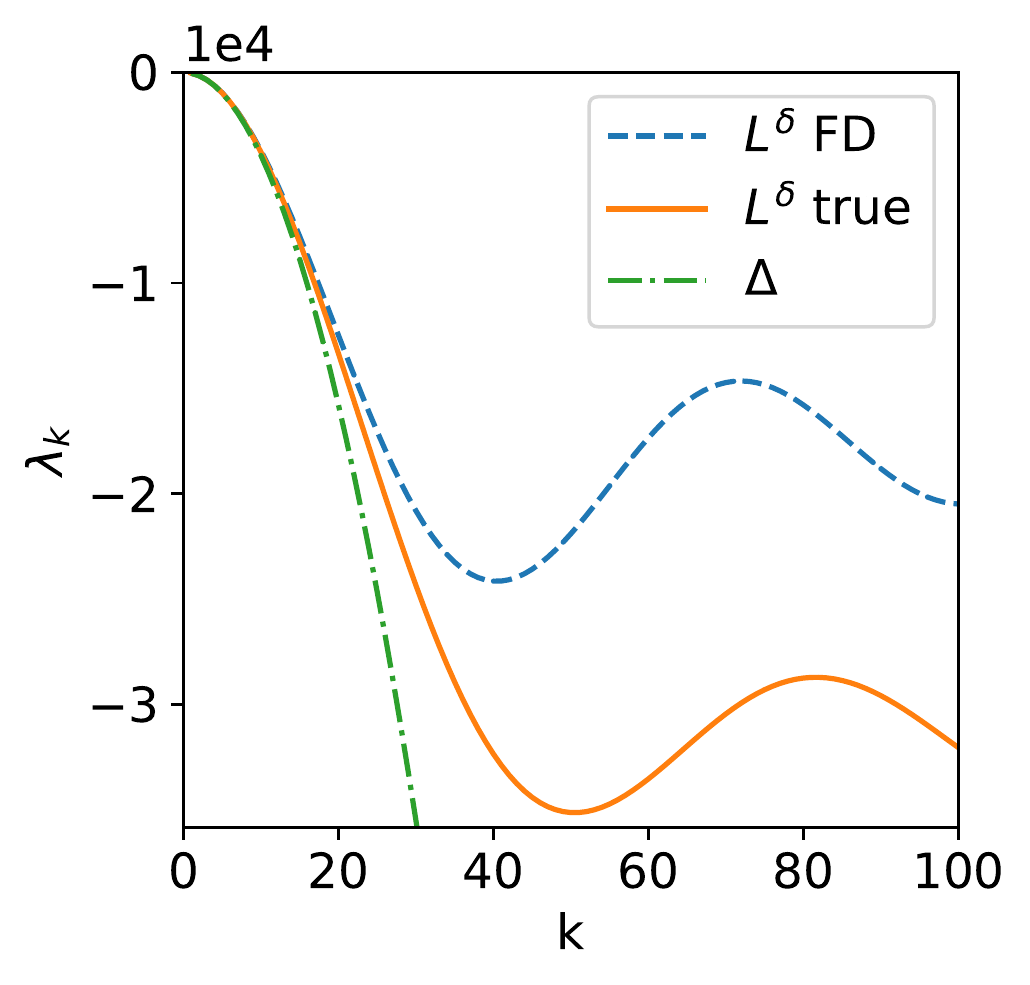}
  \end{subfigure}%
  \begin{subfigure}[t]{0.327\textwidth}
    \includegraphics[width=\textwidth]{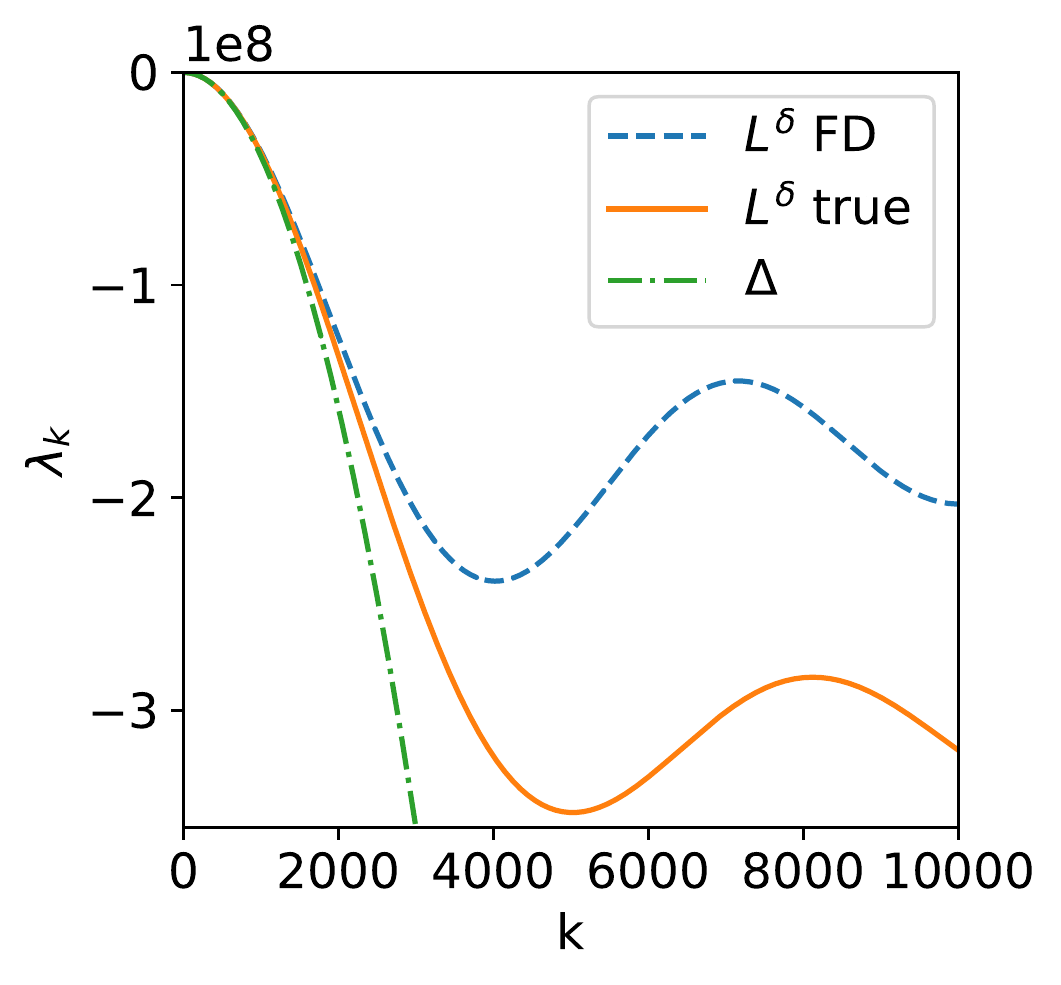}
  \end{subfigure}%
  \begin{subfigure}[t]{0.305\textwidth}
    \includegraphics[width=\textwidth]{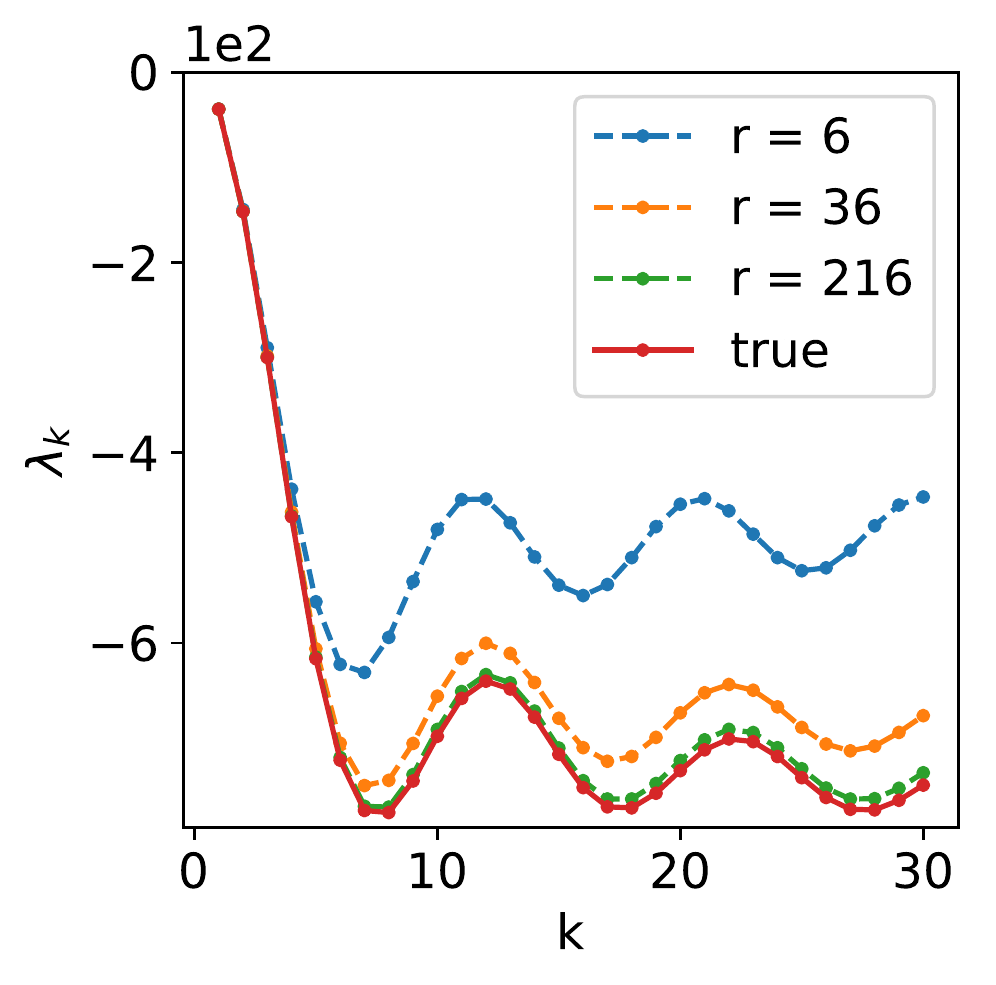}
  \end{subfigure}
  \caption{The left two figures show the eigenvalues of the finite
    difference approximation to $\Ldel$ when $\delta=3\Delta x$.
    The figures show the eigenvalues $0\le k\le N$ computed using
    $M=2N+1$ discretization points in the periodic unit interval,
    where $N=100$ and $N=10000$ respectively. The upper dashed curve
    in each figure corresponds to the eigenvalues of the finite
    difference operator.  The solid curve shows the true eigenvalue,
    using the formula above.  The lower dashed curve shows the
    corresponding eigenvalues of the 1D Laplacian for comparison.  The
    rightmost figure shows the finite difference operator's
    eigenvalues $\lambda_k$ for $1\le k\le 30$ when $L=1$ and
    $\delta=0.1$ is fixed.  The dashed curves, from top to bottom,
    correspond to discretization sizes $\Delta x = \delta/r$, for
    stencil radius $r=6$, $r=36$ and $r=216$ respectively.  The solid
    curve shows the true eigenvalues. }
  \label{fig:FD-spectrum}
\end{figure}

In this section, we compare the pseudo-spectral method in
\ref{sec:pseudo-spectral} with a finite difference approximation in 1D. Although the pseudo-spectral method described above is very accurate
for smooth, periodic problems, it does not show the same accuracy in
non-periodic situations.  For this reason,
spectral methods have typically been avoided in numerical solvers,
with finite difference and finite element methods being the preferred
treatment.  These solvers come with their own difficulties, however.
When refining the spatial mesh size $\Delta x$ to improve accuracy,
one finds that, in order for $\delta$ to remain fixed, the number of
discretization points (or finite element cells) used in computing
$\Ldel$ at a point must grow like $\delta/\Delta x$, leading to
increased computation complexity as $\Delta x$ goes to $0$.  To avoid
this blow-up in computational cost, a common technique is to fix the
ratio $\delta/\Delta x$ as $\Delta x\to 0$.  The problem with this
technique can be seen in Figure~\ref{fig:FD-spectrum}.  

In this figure, we consider the finite difference operator
$\mathbb{A}^\beta_{D,0}$ found in~\cite{tian2013analysis}, using the
scaling $\delta = r\Delta x$ with $r=3$.  This operator approximates
$\Ldel u(x)$ on an equispaced grid using a finite difference
formula of the form
\begin{equation*}
  \Ldel u(x) \approx \mathbb{A}^\beta_{D,0}u(x) = 
  a_0u(x) + \sum_{j=1}^ra_j[u(x+j\Delta x)+u(x-j\Delta x)].
\end{equation*}
When treated as a finite difference operator on a periodic interval,
the trigonometric polynomials are eigenfunctions of the discrete
operator.  To see this, let $\phi_k(x)=e^{2\pi ikx/L}$.  Then
\begin{equation*}
  \begin{split}
  \mathbb{A}^\beta_{D,0}\phi_k(x) &= \phi_k(x)\left(a_0 + 
  \sum_{j=1}^ra_j\left[e^{2\pi ikj\Delta x/L}+e^{-2\pi ikj\Delta x/L}\right]\right)\\
  &= \phi_k(x)\left(a_0+2\sum_{j=1}^ra_j\cos(2\pi kj\Delta x/L)\right).
\end{split}
\end{equation*}

The left and center plots in Figure~\ref{fig:FD-spectrum} show the
eigenvalues $\lambda_k$ (upper dashed curve) of
$\mathbb{A}^\beta_{D,0}$ associated with $1\le k\le N$ for the cases
$N=100$ and $N=10000$ respectively.  In these plots, $L=1$ and
$\Delta x = 1/(2N+1)$.  Shown for comparison are the true eigenvalues
of $\Ldel$ (solid curve) as well as the eigenvalues for the 1D
Laplacian operator (lower dashed curve).  Notice that the shape of the
eigenvalue curves do not change, indicating a certain scale
invariance.  What this indicates is that the finite difference method
with fixed $\delta/\Delta x$ ratio is \emph{never} accurate for high
frequency modes, exactly where the eigenvalues of $\Ldel$ differ
significantly from those of the 1D Laplacian.  Refining the spatial
discretization in order to better resolve high-frequency modes simply
has the effect of shifting those modes into the Laplacian-like part of
the spectrum.

The way to remedy this problem, of course, is to fix $\delta$ and
widen the finite difference stencil as $\Delta x\to 0$.  However,
truly capturing the ``interesting'' part of the spectrum typically
requires very large stencil sizes.  Consider, for example, the
right-most plot of Figure~\ref{fig:FD-spectrum}.  The figure shows the
first 30 eigenvalues of $\Ldel$ for $L=1$, $\delta=0.1$ and
$\beta=1/3$ along with the corresponding eigenvalues of the finite
difference approximations for stencils with radii $r=6$, $r=36$ and
$r=216$.  Evidently, correctly capturing the spectrum of $\Ldel$
requires finite difference stencils spanning hundreds of points.  Such
stencils are impractical in numerical solvers.

\begin{figure}
  \centering
  \includegraphics[width=0.33\textwidth]{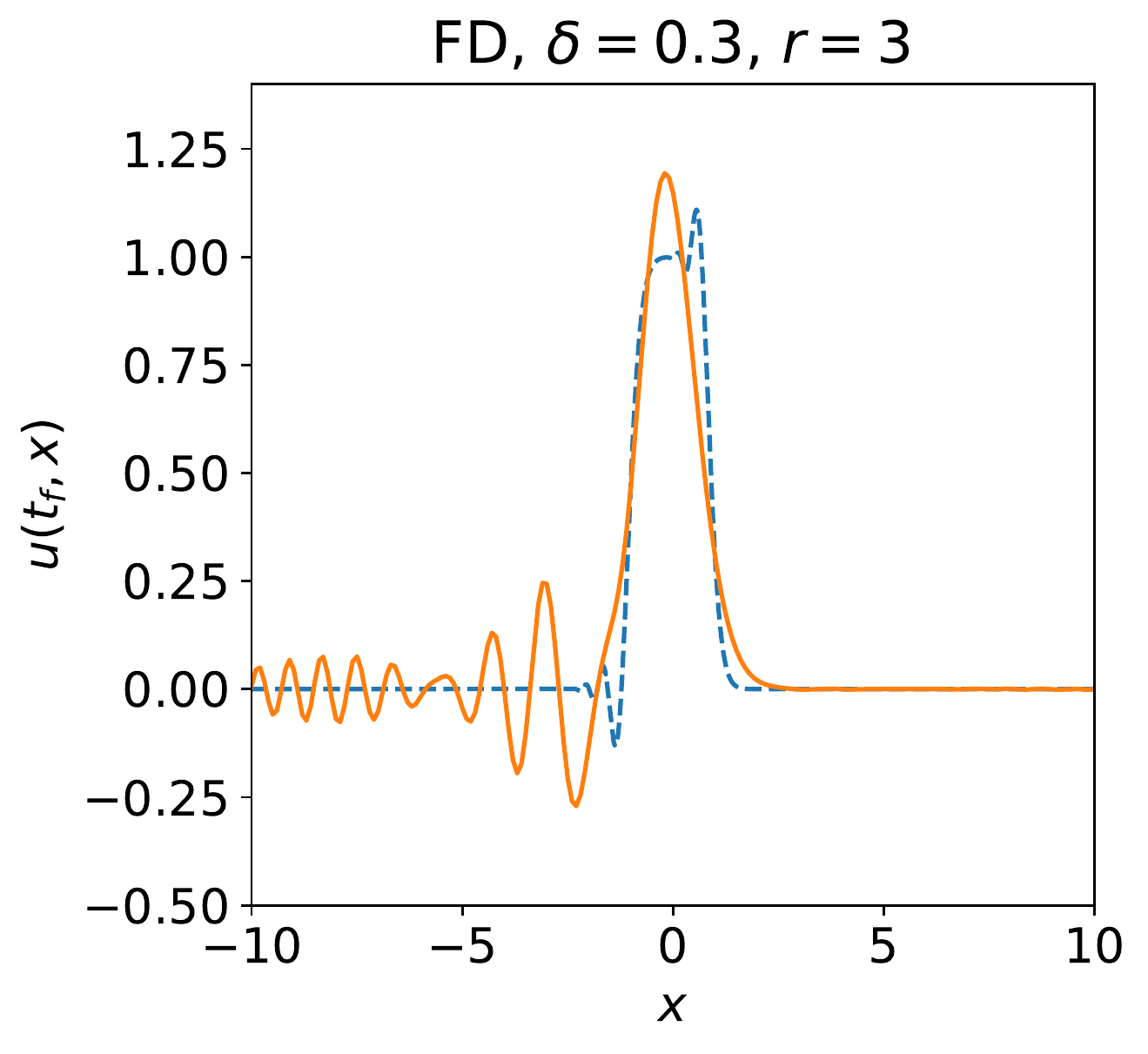}%
  \includegraphics[width=0.33\textwidth]{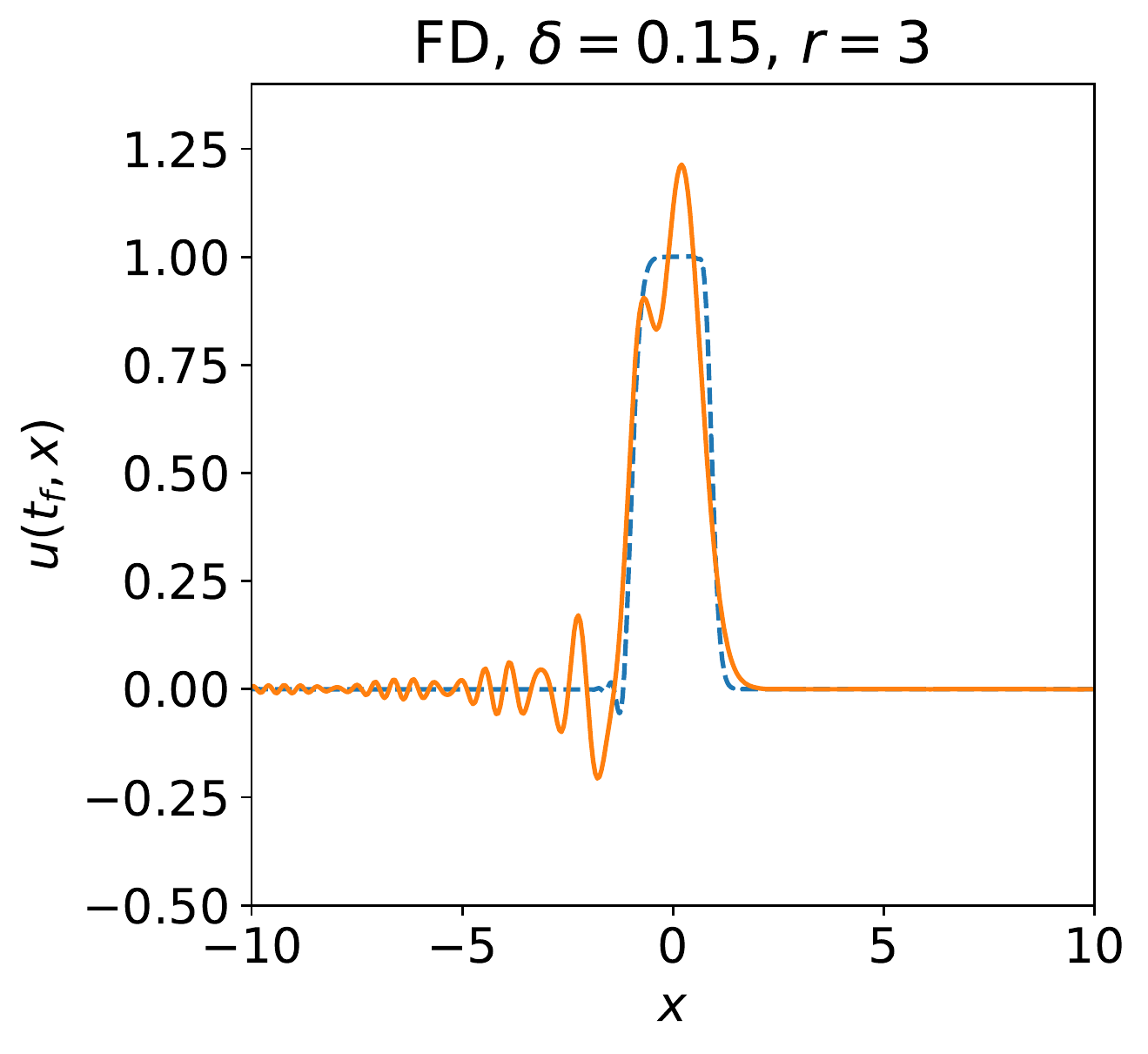}%
  \includegraphics[width=0.33\textwidth]{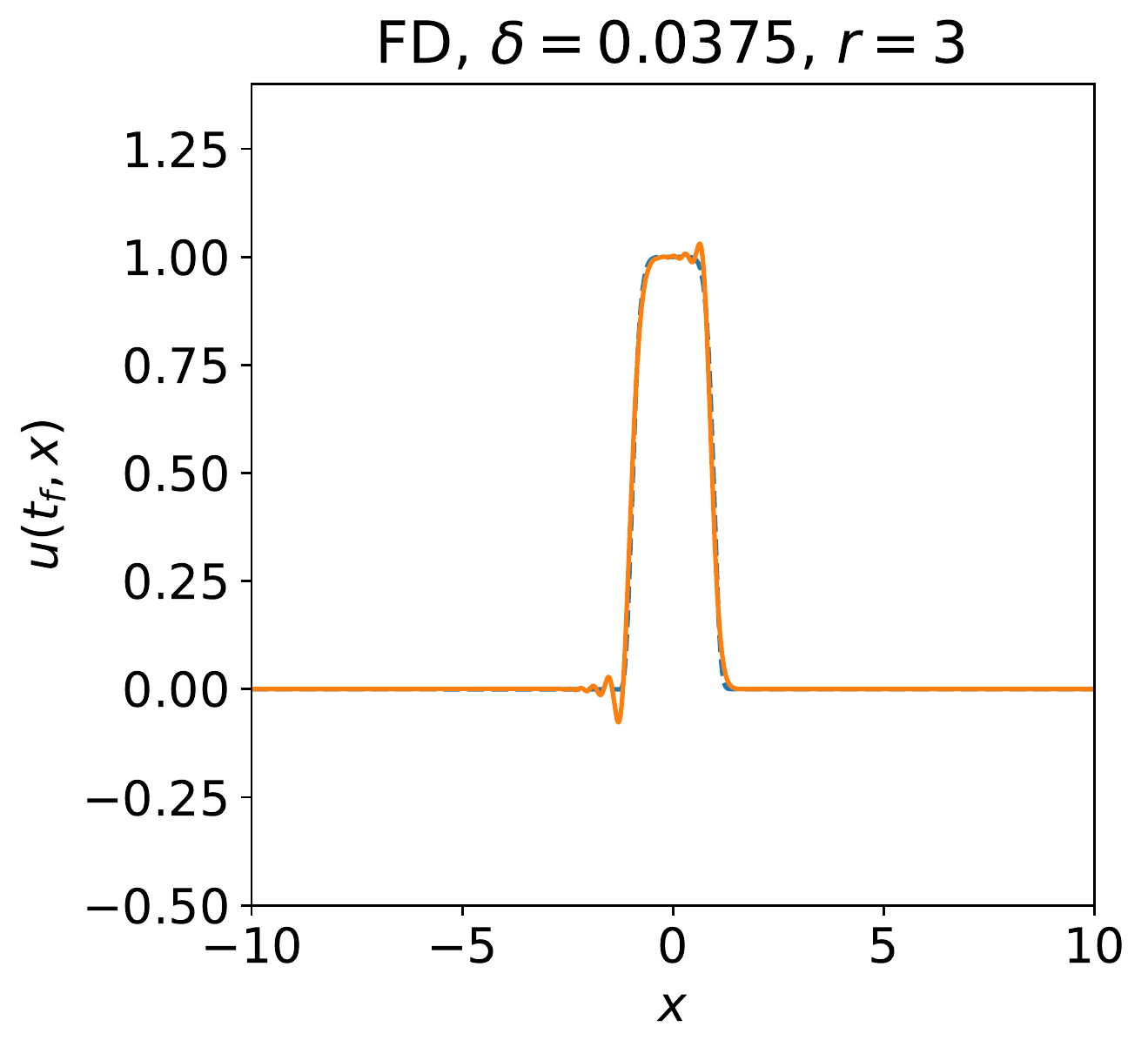}
  \includegraphics[width=0.33\textwidth]{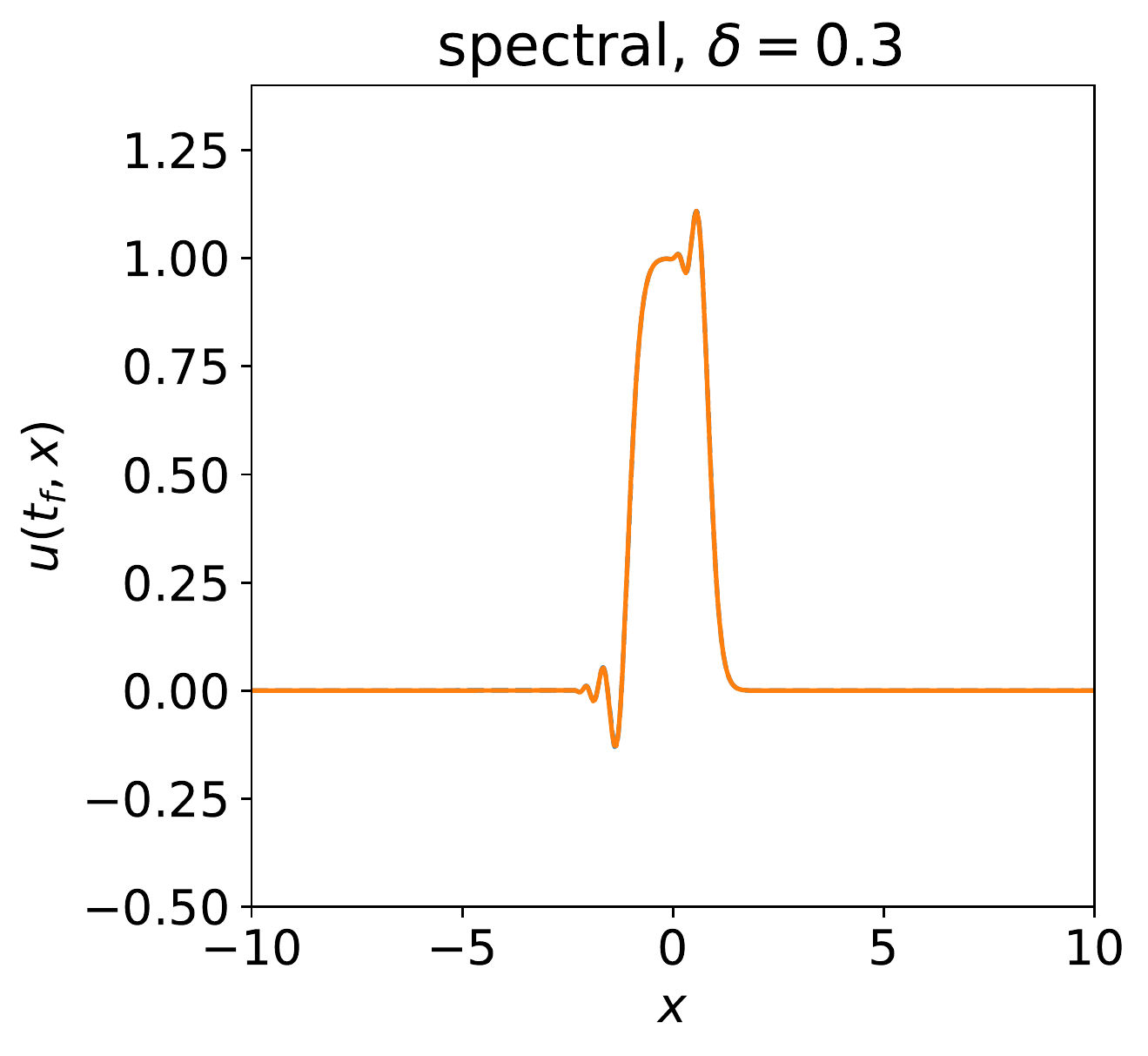}%
  \includegraphics[width=0.33\textwidth]{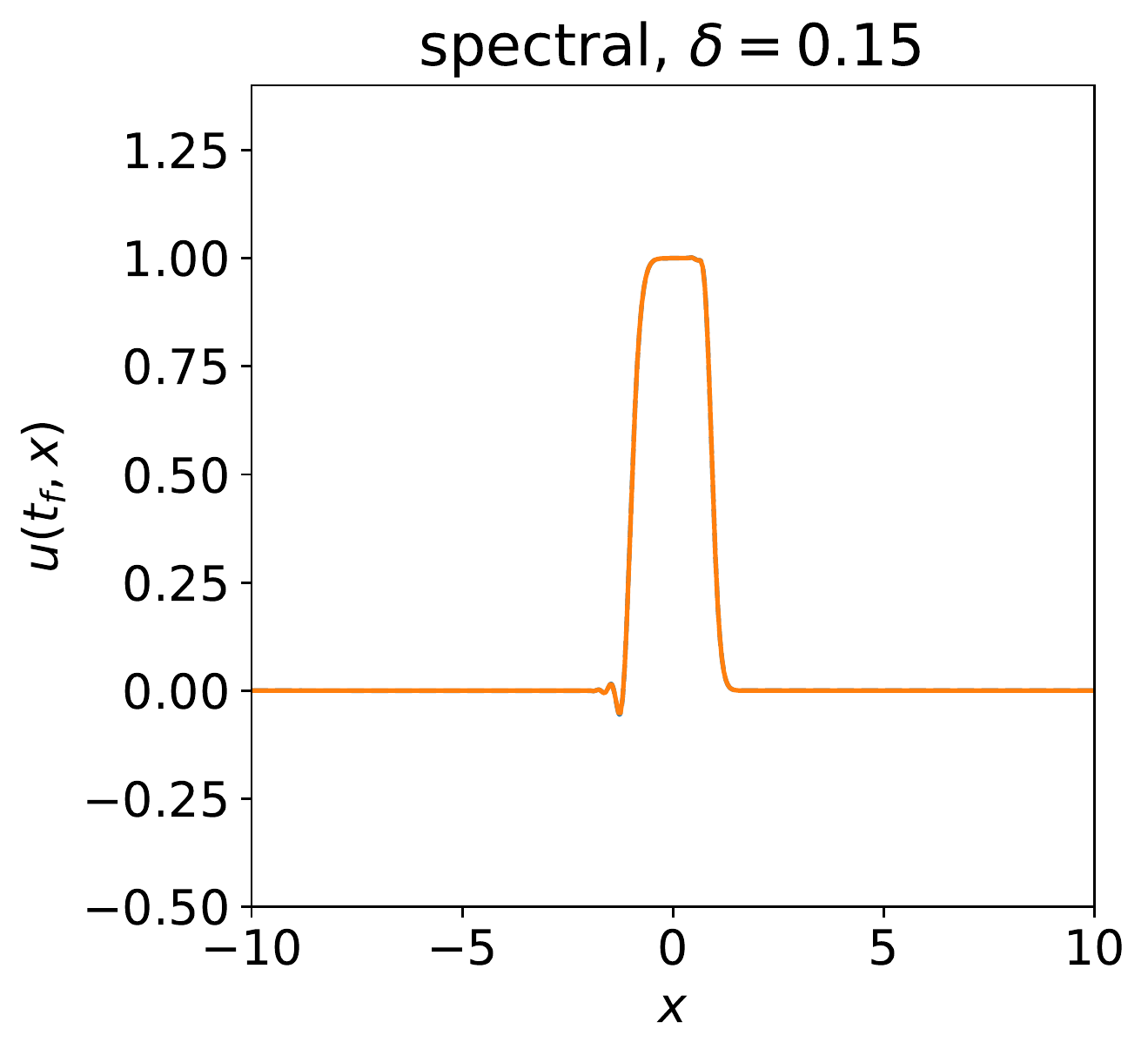}%
  \includegraphics[width=0.33\textwidth]{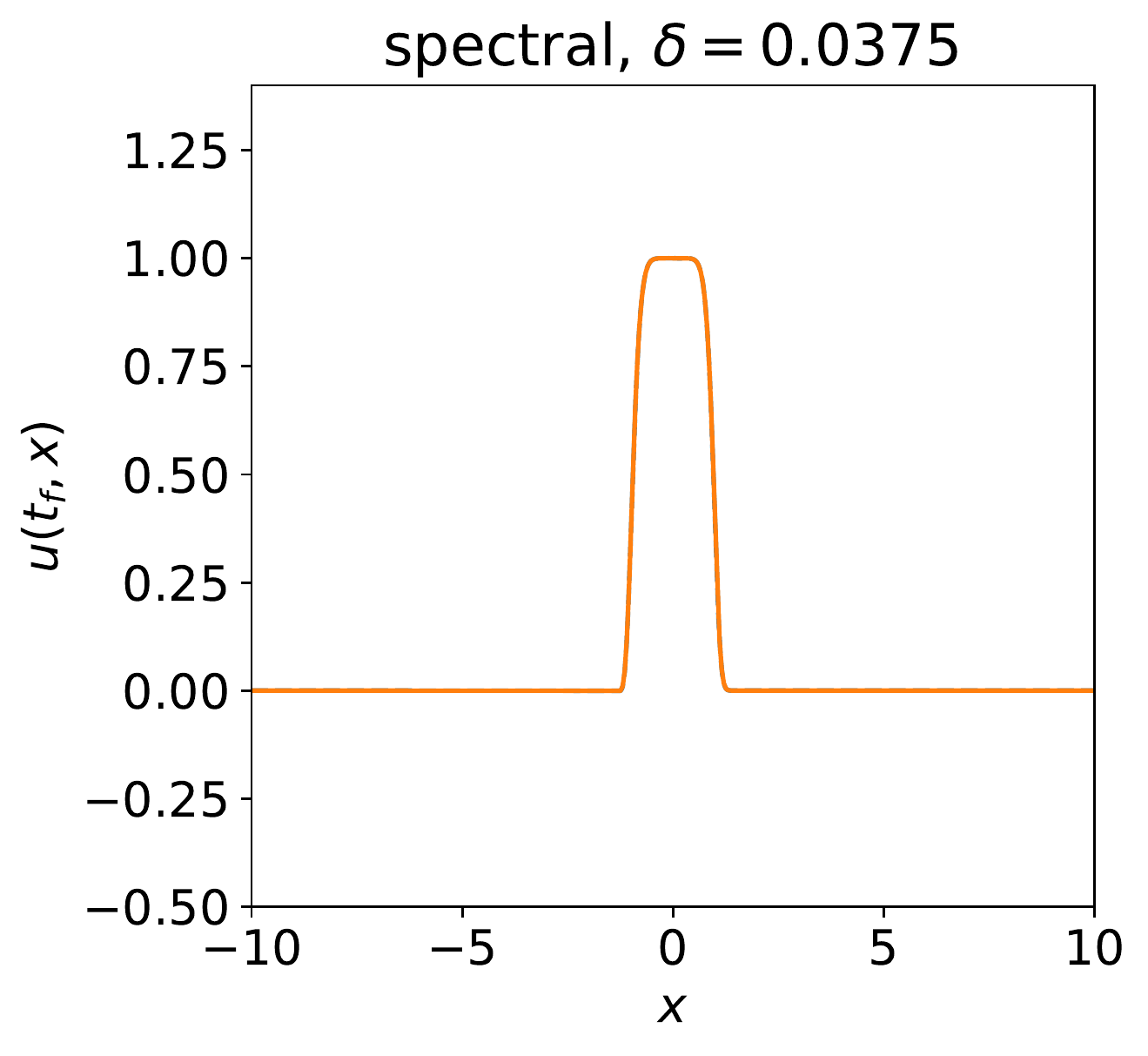}%
  \caption{Comparison of finite difference (FD) and spectral solvers
    for the nonlocal wave equation $u_{tt}=\Ldel u$.}
  \label{fig:wave-compare}
\end{figure}

\begin{figure}
  \centering
  \includegraphics[width=0.33\textwidth]{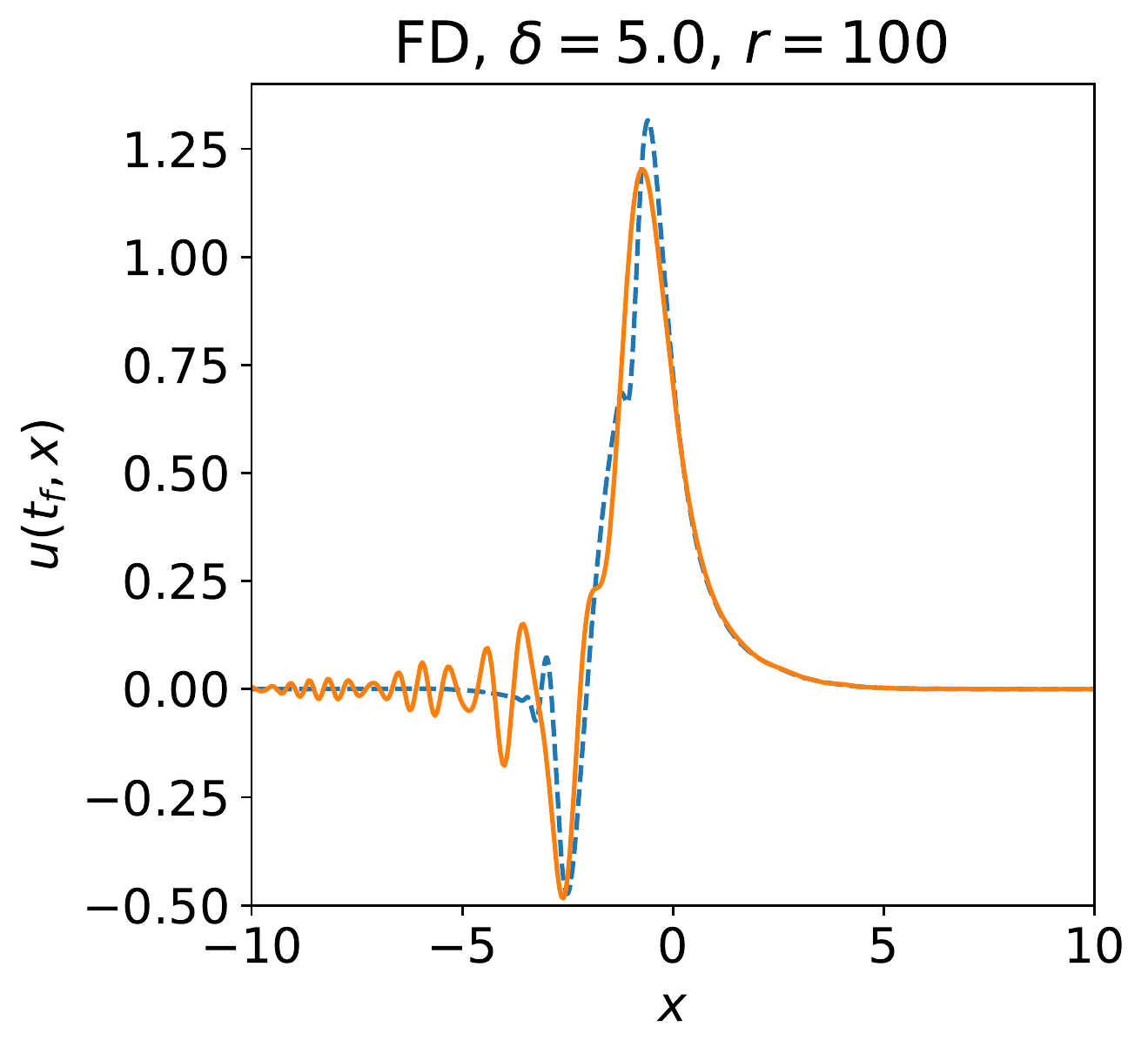}%
  \includegraphics[width=0.33\textwidth]{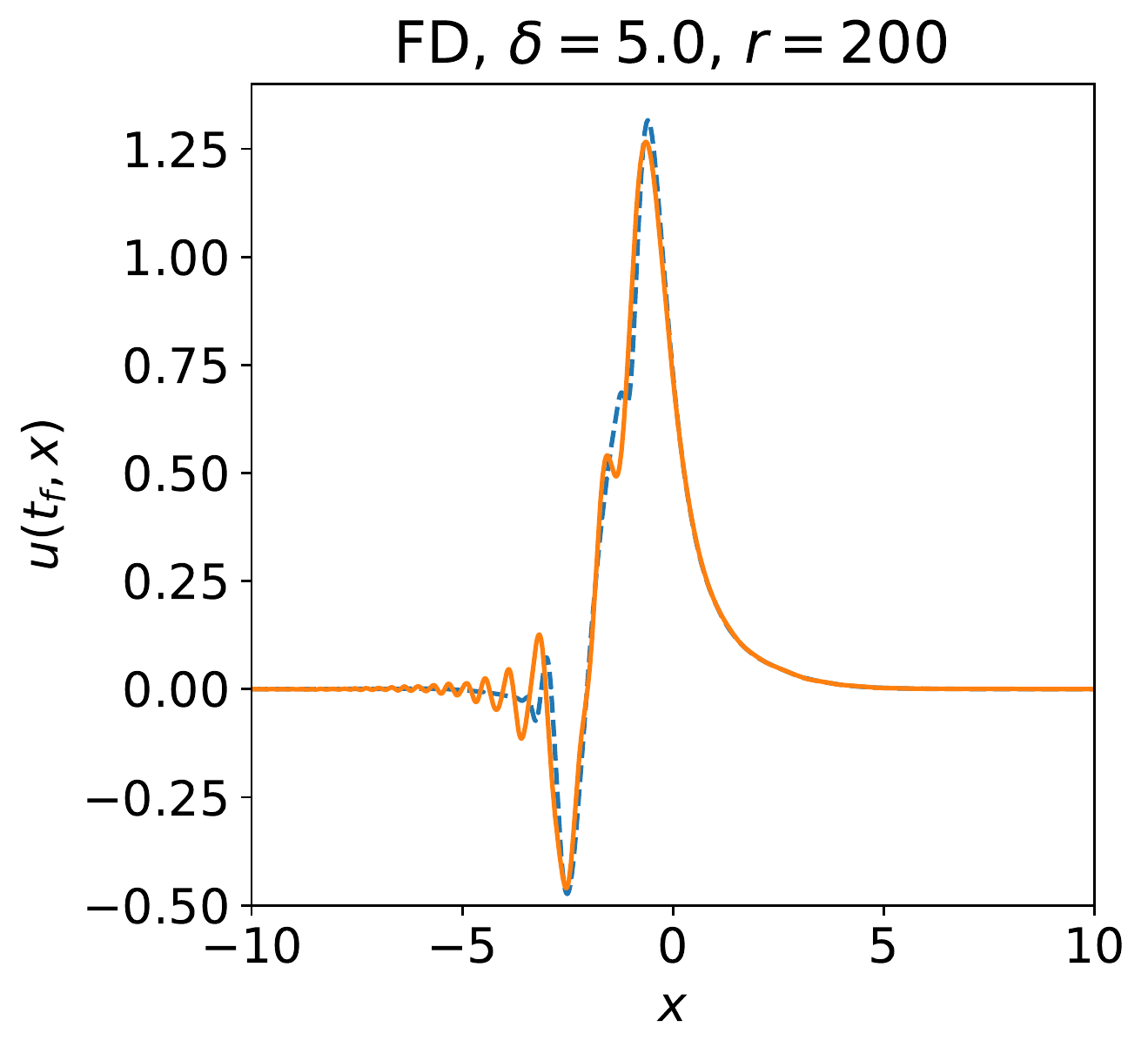}%
  \includegraphics[width=0.33\textwidth]{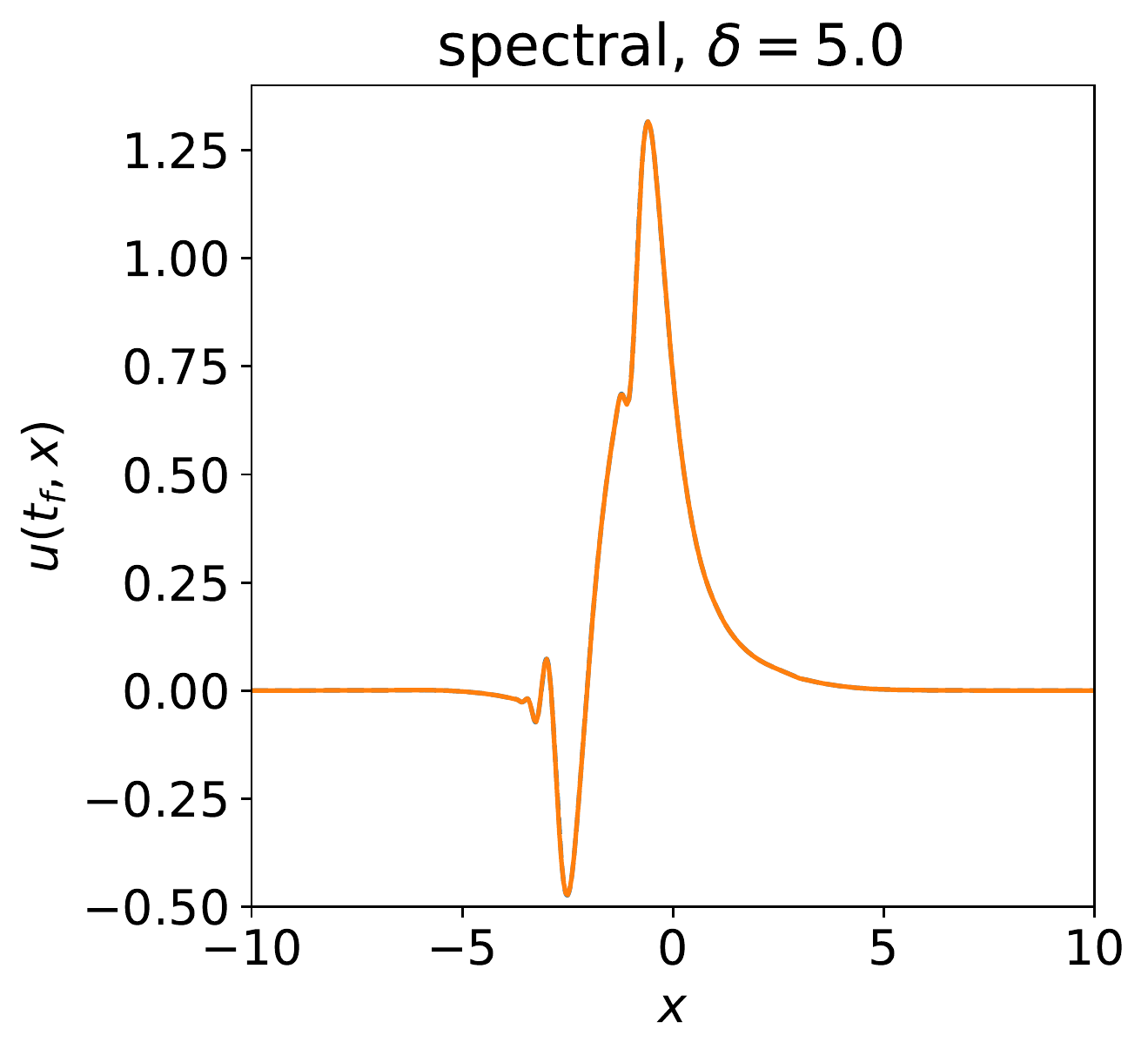}
  \caption{Spectral solutions for the nonlocal wave equation
    $u_{tt}=\Ldel u$ for large $\delta$.}
  \label{fig:wave-big-delta}
\end{figure}

As another way to see how the solutions when $\delta/\Delta x$ is
fixed can be misleading, consider Figure~\ref{fig:wave-compare}.
These plots show the results of a periodic nonlocal wave equation
\begin{gather*}
  u_{tt} = \Ldel u,\qquad (x,t)\in (0,20)\times(0,40)\\
  u(0,x) = \exp(-(x-10)^8),\quad u_t(0,x) = 0.
\end{gather*}
Note that at $x=0$ and $x=20$, the initial conditions are zero in
double precision.  For this example, rather than choosing a time stepper with fixed step size, we have instead used the \emph{solve\_ivp} method of \emph{scipy}~\cite{scipy}, which implements and adaptive RK45 timestep algorithm.

The top row shows the results of a solver based on the finite difference
$\mathbb{A}^\beta_{D,0}$ operator with $\delta=3\Delta x$,
$\beta=1/3$.
The bottom row shows the corresponding solutions for a pseudo-spectral $N=2000$ grid points.  The dashed curve in all plots shows the true solution evaluated using the FFT and~\eqref{eq:nonlocal-wave-per-sol} on a grid with $N=8000$ nodes, which is sufficient to resolve the initial data with high accuracy.

From the figure, it is clear that the
finite difference solver only approximates the true solution in the
limit as  $\Delta x$ and $\delta$ both go to zero; in other words, only when the solution
to the nonlocal equation is sufficiently close to the classical wave
equation.  In particular, this means that the standard finite
difference and finite element solvers have extreme difficulty
resolving nonlocal behaviors when the solution varies on a length
scale much smaller than $\delta$, for example, the solutions shown in
Figure~\ref{fig:wave-big-delta}, which demonstrates that, for large $\delta$ very large $r$ might be insufficient to correctly capture nonlocal wave phenomena.

\begin{table}
\centering
\begin{tabular}{|c|c|c|c|c|c|c|c|}
\hline
& \multicolumn{3}{c|}{spectral}
& \multicolumn{4}{c|}{finite differences} \\
$\delta$ & $N$ & error & CPU (s) & $N$ & $r$ & error & CPU (s) \\
\hline
0.3000 & 2000 & 7.490e-06 &   3.46 & 1000 &   3 & 4.635e-01 &   1.06 \\
0.1500 & 2000 & 6.522e-06 &   3.53 & 2000 &   3 & 3.293e-01 &   2.21 \\
0.0750 & 2000 & 7.853e-06 &   3.53 & 4000 &   3 & 1.976e-01 &   6.77 \\
0.0375 & 2000 & 9.464e-06 &   3.60 & 8000 &   3 & 9.423e-02 &  22.11 \\
\hline
5.0000 & 2000 & 7.462e-06 &   3.77 & 2000 & 100 & 3.429e-01 &  17.33 \\
& & & & 4000 & 200 & 1.934e-01 & 121.34 \\
& & & & 8000 & 400 & 9.051e-02 & 1019.27 \\
\hline
\end{tabular}
\caption{Comparison of accuracy and computing time for the spectral and finite difference methods described in Section~\ref{sec:fd-vs-spectral}.}
\label{tbl:fd-vs-spectral}
\end{table}

This observation is quantified in Table~\ref{tbl:fd-vs-spectral}, which summarizes the results displayed in Figures~\ref{fig:wave-compare} and~\ref{fig:wave-big-delta}.  For each test, the table provides the value of $N$ for the spectral method, and the value of $N$ and $r$ for the finite difference method.  Also displayed are the associated maximum absolute errors at time $t=40$ (as compared to the exact solution evaluated using the FFT and~\eqref{eq:nonlocal-wave-per-sol}) along with the total CPU time required to obtain the solution (including the time to compute the multipliers in the spectral case and time to compute the finite difference stencils in the finite difference case).  In all cases, the spectral method is able to attain an error of less than $10^{-5}$ with only $N=2000$ grid points in under 4 seconds.  On the other hand, the finite differences approach with fixed $r$ reaches an error just under $10^{-1}$ using $N=8000$ points and 22 seconds.  Worse, in the case of fixed $\delta=5$, the finite differences method with $N=8000$ points and $r=400$ obtains an error just under $10^{-1}$ and requires nearly 17 minutes.  This again shows the difficulty inherent in solving peridynamics problems accurately using finite difference or finite element approaches; one must either scale $\delta$ with the spatial discretization size (thus making the method accurate only when $\Ldel\approx\Delta$) or increase the stencil size or element order when refining the discretization, leading to very expensive computations.

\section{Discussion}
The hypergeometric representation of the Fourier multipliers \eqref{eq:multiplier-general} has shown to be a useful tool for the analysis of nonlocal equations \cite{alaliAlbin2018fourier}. 
The usefulness of this representation for computations of nonlocal equations has been demonstrated in this paper. Specifically, formula \eqref{eq:multiplier-general} allows for efficient and accurate computations of the multipliers in $n$ dimensions, without the need for integration, as shown in Section~\ref{sec:eigenvalues}. 
In addition, the fact that the nonlocality $\delta$ only appears as a parameter in the argument of the function $_2F_3$ in  \eqref{eq:multiplier-general}, leads to the uncoupling of $\delta$ and the grid size $\Delta x$. This enables the development of efficient spectral methods in which $\delta$ does not scale with $\Delta x$ as demonstrated in Section~\ref{sec:sepctral-solvers}. 
Moreover, extending the definition of the operator $\Ldel$, through extending the definition of the multipliers, to all $\beta\in\mathbb{R}\setminus\{n+4,n+6,n+8,\ldots\}$ provides a unified approach for the computations of local and nonlocal models considered in this paper.
This is demonstrated in Section~\ref{sec:sepctral-solvers}; in particular, the solutions displayed in Figures~\ref{fig:heat-example-1}, \ref{fig:heat-example-2}, (and ~\ref{fig:2d-wave-example}), which correspond to local ($\beta\ge n+2$) or nonlocal ($\beta< n+2$)
equations have been computed using the same diffusion spectral solver (wave spectral solver), respectively. Therefore, coupling local and nonlocal models is seamless using this approach.

The computational experiments in Section~\ref{sec:sepctral-solvers} indicates that, in addition to $\delta$, the parameter $\beta$ is another locality/nonlocality parameter for which smaller values of $\beta$ correspond to more nonlocal behavior while larger values correspond to more local behavior. In the context of diffusion, in Section \ref{sec:nonlocal-heat-example}, the computational experiments suggest that  $\beta$ has a critical value at $\beta=n+2$ corresponding to classical diffusion, with subdiffusion occurring when $\beta<n+2$, and superdiffusion occurring  when $\beta>n+2$. This observation combined with the smoothness of the multipliers' representation \eqref{eq:multiplier-general} in $\beta$ provide justification to extending $\Ldel$ to larger values of $\beta$.
Lastly, these observations provide a way to interpret and set apart the roles of these two nonlocality parameters with $\beta$ as the parameter that determines the type of phenomena in a given model while $\delta$ affects the scaling in the given model. 

\clearpage
\bibliographystyle{acm}
\bibliography{refs}

\end{document}